\documentclass[oneside,11pt,a4paper]{article}
\usepackage{amsfonts,epsfig,amsmath,amsthm}

\newtheorem{remark}{Remark}[section]
\newtheorem{definition}{Definition}[section]
\newtheorem{proposition}{Proposition}[section]
\newtheorem{corollary}{Corollary}[section]
\newtheorem{lemma}{Lemma}[section]
\newtheorem{theorem}{Theorem}[section]

\title{Two-scale extensions for non-periodic coefficients
\thanks{This work was carried out during the tenures of a
          fellowship from University/ITWM in Kaiserslautern (Germany)
          and an ERCIM fellowship in Luxembourg and Norway.}
}

\date{}

\author{Vsevolod Laptev\footnotemark[3]}
\begin{document}
\maketitle
\renewcommand{\thefootnote}{\fnsymbol{footnote}}
\footnotetext[3]{
NTNU-IMF,
Alfred Getz vei 1,
NO-7491 Trondheim, Norway
  (laptevv@mail.ru,\, vsevolod.laptev@math.ntnu.no)
}
\begin{abstract}
We consider non-homogeneous media with properties which can be characterized by
rapidly oscillated coefficients. 
For such coefficients we define a notion of two-scale extension,
present several ways to construct two-scale extensions,
discuss their properties and relation to homogenization
\end{abstract}


%
\section{Introduction}
It is usually difficult to predict a global behaviour of some process in 
heterogeneous media (for example composite/porous materials) although the 
physics of the process might be well understood locally. The reason lying in 
the complexity of the microstructure gives rise to different upscaling methods.

  Heterogeneities having periodic microstructure play a central role in the 
development of upscaled models. From one side they represent an important 
particular case of general heterogeneous media and on the other there are well 
developed mathematical techniques (e.g. the two-scale asymptotic expansion 
method), which help to derive formally and often rigorously the upscaled 
model. As a result many physical processes in heterogeneous media having 
periodic microstructures are well investigated both from theoretical and 
from practical points of view and the periodicity assumption is usually a 
starting point for the upscaling procedures \cite{BLP},\cite{SP},\cite{JKO}. 
Although this assumption is valid 
in only limited number of cases, mostly in artificially created materials. 
Therefore for practical purposes one should be able to deal with non-
periodic structures.

The deterministic homogenization procedure starts from 
a sequence of problems $\{{\mathcal P}^{\textstyle\varepsilon}\}$.
In the periodic case the heterogeneity in ${\mathcal P}^{\textstyle\varepsilon}$
is usually described by an $\varepsilon$-periodic function $a^\varepsilon(x)=a(x/\varepsilon)$, 
where $a(y)$ is a given $Y$-periodic function in $\mathbb{R}^d$
($Y=(0,1)^d$ is a period:
$a(y+e_i)=a(y)$, $e_i$ is a unit vector, $i=1,\dots,d$). 
Quite often the purely periodic coefficient can be generalized without 
difficulties to the locally periodic coefficient
$a^\varepsilon(x)=a(x,x/\varepsilon)$ (where $a(x,y)$ is a given 
$Y$-periodic function in $y$).
In the following steps one has to investigate the 
convergence of the sequence (in a wide sense) and to find a limit
problem $\mathcal P^0$.
The solution of the limit problem can be used in order to approximate 
the solutions of the problems ${\mathcal P}^{\textstyle\varepsilon}$ for small enough $\varepsilon$.

The coefficients $a(y)$ or $a(x,y)$ are considered in mathematical
literature as given functions belonging to some functional spaces, without
paying much attention where they come from. The construction of these
coefficients which is important for usage of homogenization 
results will be discussed in this article.

Let us assume that some process in a heterogeneous medium
occupying a bounded domain $\Omega\subset\mathbb{R}^d$ can be described by some PDE(s) with 
(at least one) rapidly oscillated coefficient $a_M(x)$, which is non necessarily periodic.
This is our initial problem ${\mathcal P}$.

Asymptotical approach applied to ${\mathcal P}$ means that we are not going to solve it
directly, but to construct a sequence of imaginary problems $\{{\mathcal P}^{\textstyle\varepsilon}\}$ 
passing through ${\mathcal P}$ at some $\bar\varepsilon$:
\begin{equation}
\begin{array}{ccccccccccc}
{\mathcal P}^{\displaystyle\varepsilon_0},&\dots&,{\mathcal P}^{\displaystyle\varepsilon_{n-1}},
&{\mathcal P}^{\displaystyle \bar\varepsilon},&{\mathcal P}^{\displaystyle\varepsilon_{n+1}},
&\dots&,{\mathcal P}^{\displaystyle \varepsilon},&
\dots&\dashrightarrow&{\mathcal P}^{0}.\\
&&&\parallel&&&&&&\\
&&&{\mathcal P}&&&&&&
\end{array}
\end{equation}
If the sequence $\{{\mathcal P}^{\textstyle\varepsilon}\}$ is convergent in some sense to 
a limit problem ${\mathcal P}^{0}$ which is easier than ${\mathcal P}$
then the solution of ${\mathcal P}^{0}$ can be used to approximate (in some sense) the solutions of 
${\mathcal P}^{\textstyle\varepsilon}$, and in particular, of 
${\mathcal P}^{\textstyle\bar\varepsilon}$ (it is our main goal). 
The ''convergence of problems'' is related to convergence 
of their solutions, but it might be restrictive  
to say something more precise.

In the periodic case, namely when $a_M(x)$ is $\bar\varepsilon$-periodic
in $\Omega$ there is a $Y$-periodic function $a(y)$  
defined in $\mathbb{R}^d$ such that $a_M(x)=a(x/\bar\varepsilon)$.
The standard sequence $\{{\mathcal P}^{\textstyle\varepsilon}\}$
is based on the $\varepsilon$-periodic coefficient $a^\varepsilon(x)=a(x/\varepsilon)$.
$a^{\bar\varepsilon}(x)=a_M(x)$ and consequently the condition 
${\mathcal P}^{\textstyle\bar\varepsilon}={\mathcal P}$
is not difficult to satisfy. This approach cannot be used for non-periodic
$a_M(x)$ since there is no such periodic $a(y)$ exists (except the case when the period 
contains the whole $\Omega$). 
But using the sequence $\{{\mathcal P}^{\textstyle\varepsilon}\}$
based on locally periodic function $a(x,y)$, where
the coefficients have the form $a^\varepsilon(x)=a(x,x/\varepsilon)$,
the requirement ${\mathcal P}^{\textstyle\bar\varepsilon}={\mathcal P}$
becomes much more realizable. We only need to find such function 
$a(x,y)$ and $\bar\varepsilon$ that $a(x,x/\bar\varepsilon)=a_M(x)$.
Therefore it is reasonable to make the following definition.
\begin{definition}
Let us say that a function $a(x,y)$, $(x,y)\in\Omega\times\mathbb{R}^d$,
$Y$-periodic in the variable $y$ 
is a two--scale extension for $a_M(x)$ if there exists a positive number
$\bar\varepsilon$ such that  
\begin{equation}
 a\left(x,\frac{x}{\bar\varepsilon}\right)
=a_M(x),\qquad \forall x\in\Omega.
\label{two-scale extension}
\end{equation}
\end{definition}
The article is organized as follows. In the next section
several ways to construct two--scale extension for arbitrary initial
coefficients $a_M(x)$ are presented. 
The Section \ref{Section Two--scale convergence and admissible test functions}  
contains a short introduction to the two-scale convergence method
together with a definition and a criterion for the concept of admissible test
function. The criterion is needed to show that the proposed in Section
\ref{Section Three approaches to construct a two--scale extension}
two scale extensions are admissible test functions in the sense of the
two-scale convergence. This is the main purpose of Sections 
\ref{Section Properties of the two-scale Continuous Extension},
\ref{Section Properties of the two-scale Discrete Extension},
\ref{s:Admissibility of the Continuous Discrete Extensions} 
(its justification consists of several results which may also be
useful of their own).
The application to the second order elliptic equation is
discussed in Section 
\ref{Section Application to the elliptic equation}. 

Why do we need this?
There are both theoretical and practical reasons to consider two-scale extensions.
First of all, they seem to be naturally related to the formal method
of two-scale asymptotic expansions and to its rigorous version --
the two-scale convergence method. If some mathematical model of a
physical process allows the formal homogenization procedure via two-scale asymptotic
expansions in the case of smooth locally periodic coefficients then as the
next step one can substitute two-scale extensions for these
coefficients and check whether the homogenization procedure remains
working for non-periodic coefficients.   

Let us now assume that our mathematical model is based on the second order
elliptic equation. The two-scale extensions might be useful for better
understanding of the following important questions related to the
concept of the averaged coefficient:\\
-- its definition, existence, properties, limits of applicability, averaging size;\\
-- connection between deterministic and stochastic approaches;\\
-- reiterative averaging (averaging of the averaged coefficient).

There are many algorithms currently known for practical calculation
of the averaged coefficient (see e.g. 
\cite{Durlofsky},\cite{Farmer},\cite{Effendiev}).
Some of them (having the same local problem with periodic boundary
conditions) can be recovered by a special choice of the two-scale extension.
This gives them a justification by an asymptotical argument as well as
some freedom for improvement and generalization. For example it is
possible to correct the averaged solution in a postprocessing step
using a standard technique from homogenization theory \cite[p.76]{BLP}.
Therefore for the practical problems like heat transfer in composite
materials and unsaturated flow in heterogeneous porous media the  
choice of the two-scale extension defines a numerical method which can
be used as a possible alternative to such methods as multiscale finite element method 
\cite{Hou},\cite{Cai} or heterogeneous multiscale method \cite{EE}.
\section{Three approaches to construct a two--scale extension}
\label{Section Three approaches to construct a two--scale extension}
First of all we have the ${\mathcal Trivial}$ Extension:
$$
  a(x,y):=a_M(x),\qquad x\in\Omega, y\in\mathbb{R}^d.
$$
But we cannot expect something better than the constant sequence
$\{{\mathcal P}^{\textstyle\varepsilon}\}=\{{\mathcal P}\}$ with the limit problem 
${\mathcal P}^{0}={\mathcal P}$ which is just as difficult to solve.
This practically useless extension gives although an approximation to ${\mathcal P}$
with a perfect quality. Different two-scale extensions lead to
upscaled problems with different quality. At least we know that not
all are bad.

For the other two approaches we need to know $a_M(x)$ in a neighbourhood
of a point in $\Omega$.
Since this can create some problems close to the boundary,
let us assume that $a_M(x)$ can be somehow extended to 
a larger domain $\tilde\Omega$ which is also bounded (if we find nothing better, we can
choose some value of $a_M(\cdot)$ in $\Omega$ as a constant value in
$\tilde\Omega\setminus\Omega$).

Next we need to choose $\bar\varepsilon$. For periodic $a_M(x)$ it is 
reasonable to choose $\bar\varepsilon$ equal to the period, but in general
we are free in choosing it. 
Let $W(x)$ be an $\bar\varepsilon$--cube with the center $x$ and sides aligned with the coordinate axes.
Up to now the only restrictions on $\bar\varepsilon$ are:
we consider $\bar\varepsilon$ to be small comparing to the typical size of $\Omega$
and all cubes $W(x)$, $x\in\Omega$ should be completely inside $\tilde\Omega$. 

Having in mind the volume averaging method it might be reasonable to call 
$W(x)$ as a (cubic) representative elementary volume (REV) around the point $x$. 

Two approaches to construct the two--scale extension $a(x,y)$ for $a_M(x)$
are different in the sense that the first is created via continuous 
(${\mathcal Continuous}$ Extension) and the second via discrete 
(${\mathcal Discrete}$ Extension) 'motion' of $W(x)$ in $\Omega$.
\subsection{${\mathcal Continuous}$ Extension}
\label{Section Continuous Extension}  
Let $x$ be some fixed point in $\Omega$.\\
$\bullet$ First we define an auxiliary function $\tilde a(x,\cdot)$ at $y\in W(x)$:
$$
\tilde a(x,y)=a_M(y),\qquad y\in W(x).
$$
$\bullet$ Secondly we extend it to the whole $\mathbb{R}^d$ periodically --
$\tilde a(x,y)$ is $\bar\varepsilon$-periodic in $y$.\\ 
$\bullet$ Thirdly
$$
a(x,y):=\tilde a(x,\bar\varepsilon y)
$$
is defined in $\Omega\times\mathbb{R}^d$, $Y$--periodic in $y$. 
It satisfies (\ref{two-scale extension}).

\subsection{${\mathcal Discrete}$ Extension}
\label{Section Discrete Extension}
Let us assume that we have some finite partition 
$\overline\Omega=\cup_j\overline\Omega_j$, $j=1,\dots ,N$.
$\Omega_i\cap\Omega_j=\emptyset$, $i\ne j$. 
For each $\Omega_j$ there is a corresponding $\bar\varepsilon$--cube 
$W_j=W(\hat x^j)$, $\Omega_j\subseteq W_j$, $\hat x^j$ is a center of $W_j$.
\\ 
$\bullet$ First for any fixed $x\in\Omega\cap\Omega_j$ we define an auxiliary 
function $\tilde a(x,\cdot)$ at $y\in W_j$:
$$
  \tilde a(x,y):=a_M(y),\qquad y\in W_j.
$$
$\bullet$ Secondly we extend it to the whole $\mathbb{R}^d$ periodically --
$\tilde a(x,y)$ is $\bar\varepsilon$--periodic in $y$.\\
$\bullet$ Thirdly
$$
a(x,y):=\tilde a(x,\bar\varepsilon y)
$$
is defined in $\Omega\times\mathbb{R}^d$, $Y$--periodic in $y$.
It satisfies (\ref{two-scale extension}).

\begin{remark} 
Both extensions are also well defined in $\overline\Omega\times Y$ 
(this will help to show continuity of some properties in 
$\overline\Omega$). For
$\Omega_j$ much smaller than $W_j$, $\hat x^j\in\Omega_j$ 
the ${\mathcal Discrete}$ Extension
can be seen as a discretization of the ${\mathcal Continuous}$ Extension. 
\end{remark}
In order to use the results of convergence and error estimations,
one usually needs smoothness of $a(x,y)$. However it is easy to see that 
$a(x,y)$ are continuous neither in $x$ nor in $y$ 
(and are properly defined only a.e.).
Anyway, in the next sections our goal will be to show that these
$a(x,y)$ can be considered as admissible test functions in the sense 
of two--scale convergence and at least for the second order elliptic 
equation with highly oscillated (conductivity, permeability) coefficient the standard 
procedure ~\cite{Allaire} still works and solutions of 
$\{{\mathcal P}^{\textstyle\varepsilon}\}$ converge to the solution 
of ${\mathcal P}^{0}$. 

Please note that the convergence of the solutions of
$\{{\mathcal P}^{\textstyle\varepsilon}\}$
is important, but it cannot guarantee that the solution of ${\mathcal P}$
can be well--approximated with the help of the solution of the problem
${\mathcal P}^{0}$.
The approximation may fail since ${\mathcal P}^{\textstyle\bar\varepsilon}$ plays 
a central role in the construction of the sequence and even if the sequence 
"converges", ${\mathcal P}^{0}$ may be 'close' to practically useless 
problems ${\mathcal P}^{\textstyle\varepsilon}$, for 
$\varepsilon\ll\bar\varepsilon$ but still 'far' from 
${\mathcal P}^{\textstyle\bar\varepsilon}$.

\subsection{An example of ${\mathcal P}^{\textstyle\varepsilon}$ for the 
elliptic problem ${\mathcal P}$}
\label{Subsection Elliptic example}

In this example we consider the second order elliptic problem with homogeneous
Dirichlet boundary condition as the initial problem
\begin{equation}
\label{initial elliptic problem}
\mbox{${\mathcal P}$:}\qquad
-\nabla\cdot(a_M(x)\nabla u)=f\qquad\mbox{in $\Omega$,}
\qquad u|_{\partial\Omega}=0;
\end{equation}
and the sequence of problems $\{{\mathcal P}^{\textstyle\varepsilon}\}$ is
\begin{equation}
\mbox{${\mathcal P}^{\textstyle\varepsilon}$}:\qquad
-\nabla\cdot(a(x,x/\varepsilon)\nabla u_\varepsilon)=
f\qquad\mbox{in $\Omega$,}
\qquad u_\varepsilon|_{\partial\Omega}=0,
\label{sequence P(varepsilon)}
\end{equation}
where 
$a_M(x)=\{a_M^{ij}(x)\}$ and $a(x,y)=\{a^{ij}(x,y)\}$ are $d\times d$
matrix functions in general case and $a^{ij}(x,y)$ is a two--scale extension of $a_M^{ij}(x)$.

Naturally $a_M(\cdot)$ is required to be bounded and positive definite.
Can we expect similar properties for $a(x,x/\varepsilon)$ which are 
important for verification that $\{{\mathcal P}^{\textstyle\varepsilon}\}$
is a sequence of solvable problems?
 
\subsection{Properties of ${\mathcal Continuous}$ and  
${\mathcal Discrete}$ Extensions inherited from $a_M(x)$}
\begin{proposition}
A property of $a_M(x)$ which is valid for all $x\in\tilde\Omega$
is also valid for $a(x,y)$ in $\Omega\times Y$. 
\end{proposition}
\begin{proof}
For both ${\mathcal Continuous}$ Extension and 
${\mathcal Discrete}$ Extension there is a mapping $z: \Omega\times
Y\longrightarrow\tilde\Omega$ that $a(x,y)=a_M(z(x,y))$. 
\end{proof}
\begin{corollary}
\label{P0} 
Let ${\mathcal M}$ be the mapping ${\mathcal M}: a_M(\cdot)\longrightarrow a(\cdot,\cdot)$.
Then\\
$\bullet$ ${\mathcal M}$ is linear.\\
$\bullet$ $|a_M(\cdot)|\stackrel{\mathcal M}{\longrightarrow} |a(\cdot,\cdot)|$.\\
$\bullet$ $a_M(\cdot)^p\stackrel{\mathcal M}{\longrightarrow} a(\cdot,\cdot)^p$.\\
$\bullet$ if $a_M(x)$ is uniformly bounded, positive definite matrix function in 
$\tilde\Omega$, $a_M^{ij}(\cdot)\stackrel{\mathcal M}{\longrightarrow}a^{ij}(\cdot,\cdot)$
then $a(x,y)$ is uniformly bounded, positive definite matrix
function in $\Omega\times Y$.
\end{corollary}
\begin{proof}
For example, if $b_M(\cdot)=|a_M(\cdot)|$ then\\ 
{$b(x,y)=b_M(z(x,y))=|a_M(z(x,y))|=|a(x,y)|$}. Similar with others.
\end{proof}
We note that ${\mathcal M}$ for the ${\mathcal Discrete}$ Extension
has some similarity with the unfolding operator ${\mathcal T}$ \cite{Cioranescu}.
\section{Two--scale convergence and admissible test functions}
\label{Section Two--scale convergence and admissible test functions}
The concept of two--scale convergence was introduced in ~\cite{Ng} and
further developed in ~\cite{Allaire}. A recent review of a two--scale
convergence in $L^p(\Omega)$ space can be found in ~\cite{Lukkassen}.
In this section we formulate some results
related to two--scale convergence in $L^2(\Omega)$ mainly following ~\cite{Allaire},
but with some modifications of the concept of admissible test function.
We will need these results in 
Section~\ref{s:Admissibility of the Continuous Discrete Extensions}.

\begin{definition}
Let ${\cal B}_{TF}={\cal D}(\Omega\times Y)$ be a base space of test functions.
\end{definition}
A function $f(x,y)$ initially defined a.e. in 
$\Omega\times\overline Y$ we can extend to
a $Y$--periodic function in $\Omega\times \mathbb{R}^d$ by periodical repetition,
except perhaps the points periodic to $\partial Y$.
\begin{lemma}
\label{lem Continuous}
For any $Y$--periodic function $\psi(x,y)\in C(\overline\Omega\times\overline Y)$
\begin{equation}
\lim\limits_{\varepsilon\to 0}\int_{\Omega}\psi(x,\frac{x}{\varepsilon})\,dx=
\int_{\Omega}\int_{Y}\psi(x,y)\,dx\,dy
\end{equation}
\end{lemma}
\begin{proof}
For example see~\cite{Radu}.
\end{proof}

In the following we will deal with sequences $\{u_\varepsilon\}$.
$u_\varepsilon$ is a pair $(u,\varepsilon)\in L^2(\Omega)\times\mathbb{R}_+$.
The sequence $\{u_\varepsilon\}$ is a sequence of pairs $\{(u_n,\varepsilon_n)\}_{n=0}^\infty$
where $\{\varepsilon_n\}$ is a fixed sequence of strictly positive numbers tending to zero.
"$\lim\limits_{\varepsilon\to 0}$" is the same as 
"$\lim\limits_{\scriptstyle n\to\infty\atop\scriptstyle\varepsilon=\varepsilon_n}$".
\begin{definition}
\label{TSC}
A sequence $\{u_\varepsilon(x)\}$ from $L^2(\Omega)$ is said to be two-scale convergent to a
limit $u_0(x,y)\in L^2(\Omega\times Y)$ if 
\begin{itemize}
\item[(i)] \label{TSC1} for all $\psi\in{\cal B}_{TF}$:
\begin{equation}
\lim\limits_{\varepsilon\to 0}\int_{\Omega}
u_\varepsilon(x)\psi(x,\frac{x}{\varepsilon})\,dx=
\int_{\Omega}\int_{Y}u_0(x,y)\psi(x,y)\,dx\,dy
\label{TSCeq}
\end{equation}
\item[(ii)] \label{TSC2} $u_\varepsilon$ is bounded in $L^2(\Omega).$
\end{itemize}
\end{definition}
We prefer to insure that all two-scale convergent
sequences are bounded. Having chosen ${\cal B}_{TF}$ somewhat larger, for instance 
$L^2[\Omega;C_{per}(Y)]$ we would have (i) $\Rightarrow$ (ii) due to weak convergence 
of $u_\varepsilon$. We refer to ~\cite{Lukkassen} for the discussion
of this topic and for the definitions of the functional spaces like $L^2[\Omega;C_{per}(Y)]$. 
\begin{remark}
\label{TSC has sence}
The Def.~\ref{TSC} has sense since the limit $u_0(x,y)$ is unique as an element of 
$L^2(\Omega\times Y)$ 
due to density of ${\cal B}_{TF}$ in 
$L^2(\Omega\times Y)$ and at least the following 
sequences are two-scale convergent:
\begin{enumerate}
\item If $\phi(x,y)\in C(\overline\Omega\times\overline Y)$ then 
$u_\varepsilon(x)=\phi(x,x/\varepsilon)$ two-scale converges to $\phi(x,y)$.
\item If $u_\varepsilon(x)\to u(x)$ in $L^2(\Omega)$ then $u_\varepsilon(x)$
two-scale converges to $u_0(x,y)=u(x)$.
\end{enumerate}
\end{remark}
\begin{proof}
In both cases $u_\varepsilon$ is bounded. 
The first statement is a consequence of Lem.~\ref{lem Continuous}.
To use Lem.~\ref{lem Continuous} in the second statement we should
approximate $u(x)$ by a smooth function in $L^2(\Omega)$. 
\qquad
\end{proof}
\begin{remark}
Usually in the definition of the two--scale convergence one uses 
${\cal B}_{TF}={\cal D}[\Omega;C^\infty_{per}(Y)]$.
\end{remark}

If we want to check that some sequence $\{u_\varepsilon\}$ is two--scale convergent
then it is better to have possibly smaller set of test functions 
(${\cal B}_{TF}$). But if we already know that $\{u_\varepsilon\}$ is two--scale
convergent (for example from compactness result, see Cor.~\ref{compactness})
then it is desirable to be much more free in choosing $\psi$ for
(\ref{TSCeq}).
\begin{definition}
\label{ATF}
A $Y$--periodic function $\phi(x,y)$ 
square integrable in $\Omega\times Y$ 
with well defined $\phi(x,x/\varepsilon)$ in $L^2(\Omega)$ 
for all $\varepsilon\in\{\varepsilon_n\}$
is called an admissible test function (ATF)
if for all two-scale convergent sequences 
$\{u_\varepsilon\}$ with a limit $u_0(x,y)$ holds:
$$
\lim\limits_{\varepsilon\to 0}\int_{\Omega}
u_\varepsilon(x)\phi(x,\frac{x}{\varepsilon})\,dx=
\int_{\Omega}\int_{Y}u_0(x,y)\phi(x,y)\,dx\,dy.
$$
\end{definition}
Please note that we do not consider $\phi$ to be an element of
$L^2(\Omega\times Y)$ since different representatives $\bar\phi(x,y)$, $\bar{\bar\phi}(x,y)$
of the same element $\phi\in L^2(\Omega\times Y)$ may have
$\bar\phi(x,\frac{x}{\varepsilon})\ne\bar{\bar\phi}(x,\frac{x}{\varepsilon})$
in $L^2(\Omega)$, or even $\bar\phi(x,\frac{x}{\varepsilon})\notin L^2(\Omega)$. 
\begin{theorem}(see ~\cite{Allaire}, Th.1.2.)
\label{compactness, Allaire}
From any bounded sequence $\{u_\varepsilon\}$ in $L^2(\Omega)$ it is
possible to extract a subsequence $\{u_\varepsilon'\}$ and there exists 
$u_0(x,y)\in L^2(\Omega\times Y)$ so that for all 
$\phi(x,y)\in L^2[\Omega;C_{per}(Y)]$:
$$
\lim\limits_{\varepsilon\to 0}\int_{\Omega}
u_\varepsilon'(x)\phi(x,\frac{x}{\varepsilon})\,dx=
\int_{\Omega}\int_{Y}u_0(x,y)\phi(x,y)\,dx\,dy.
$$
\end{theorem}
\begin{corollary}
\label{compactness}
From any bounded sequence $\{u_\varepsilon\}$ in $L^2(\Omega)$ it is
possible to extract a two-scale convergent subsequence. 
\end{corollary}
\begin{proof}
${\mathcal B}_{TF}=D(\Omega\times Y)\subset L^2[\Omega;C_{per}(Y)]$.\qquad
\end{proof}
\begin{corollary}
All functions from $L^2[\Omega;C_{per}(Y)]$ are ATF.
\end{corollary}
\begin{proof}
Let us assume the opposite: $\phi(x,y)\in L^2[\Omega;C_{per}(Y)]$ and $u_\varepsilon(x)$
two-scale converges to $u_0(x,y)$, but there exists $\delta>0$, subsequence 
$\{u_\varepsilon'\}$ that
\begin{equation}
\label{contr assumption}
|\int_\Omega u_\varepsilon'(x)\phi(x,x/\varepsilon)\,dx-
\int_\Omega\int_Y u_0(x,y)\phi(x,y)\,dx\,dy|\ge\delta.
\end{equation}
From Def.~\ref{TSC}(ii), Th.~\ref{compactness, Allaire}
there exists a subsequence $u_\varepsilon''$ in $u_\varepsilon'$ that
for all $\psi(x,y)\in L^2[\Omega;C_{per}(Y)]$
$$
\lim\limits_{\varepsilon\to 0}\int_{\Omega}
u_\varepsilon''(x)\psi(x,\frac{x}{\varepsilon})\,dx=
\int_{\Omega}\int_{Y}u_1(x,y)\psi(x,y)\,dx\,dy.
$$
$u_1(x,y)=u_0(x,y)$ due to the uniqueness of the two-scale limit of 
$\{u_\varepsilon''\}$. For $\psi=\phi$ there is a contradiction with 
(\ref{contr assumption}).\qquad
\end{proof}
\subsection{Necessary and sufficient conditions for $\phi$ to be ATF}
Let us assume that $\phi(x,y)\in{\cal A}_{TF}$ -- a set of ATF. 
\begin{itemize}
\item
First we test Def.~\ref{ATF}
with $u_\varepsilon(x)=\psi(x,x/\varepsilon)$, for all $\psi\in{\cal B}_{TF}$.
$\{u_\varepsilon\}$ two scale converges to $\psi(x,y)$ 
(see Rem.~\ref{TSC has sence})
$$
\lim\limits_{\varepsilon\to 0}\int_{\Omega}
\psi(x,\frac{x}{\varepsilon})\phi(x,\frac{x}{\varepsilon})\,dx=
\int_{\Omega}\int_{Y}\psi(x,y)\phi(x,y)\,dx\,dy.
$$
\item
Second we test Def.~\ref{ATF} with $u_\varepsilon(x)=u(x)$, for all 
$u(x)\in L^2(\Omega)$ (see Rem.~\ref{TSC has sence}). 
$$
\lim\limits_{\varepsilon\to 0}\int_{\Omega}
u(x)\phi(x,\frac{x}{\varepsilon})\,dx=
\int_{\Omega}\int_{Y}u(x)\phi(x,y)\,dx\,dy.
$$
then $\phi(x,x/\varepsilon)$ weakly converges to $\int_Y\phi(x,y)\,dy$ in
$L^2(\Omega)$ and consequently is bounded.
From 1 and 2 we conclude that $u_\varepsilon(x)=\phi(x,x/\varepsilon)$
two-scale converges to $\phi(x,y)$.
\item
Third we test Def.~\ref{ATF} with $u_\varepsilon(x)=\phi(x,x/\varepsilon)$:
$$
\lim\limits_{\varepsilon\to 0}\int_{\Omega}
\phi(x,\frac{x}{\varepsilon})^2\,dx=
\int_{\Omega}\int_{Y}\phi(x,y)^2\,dx\,dy.
$$
\end{itemize}
\begin{proposition}
\label{NC} 
The necessary conditions for a function $\phi$ to be from ${\cal A}_{TF}$:
\begin{subequations}
\label{ATF conditions}
\begin{equation}
\mbox{$\phi(x,x/\varepsilon)$ two-scale converges to $\phi(x,y)$}
\label{ATF condition 1}
\end{equation}
\begin{equation}
\mbox{$\lim\limits_{\varepsilon\to 0}\|\phi(x,x/\varepsilon)\|_{L^2(\Omega)}=
\|\phi(x,y)\|_{L^2(\Omega\times Y)}$}
\label{ATF condition 2}
\end{equation}
\end{subequations} 
\end{proposition}
The conditions implicitly require that $\phi(x,y)$ is square integrable in 
$\Omega\times Y$ and $\phi(x,x/\varepsilon)\in L^2(\Omega)$ is 
well-defined for all $\varepsilon\in\{\varepsilon_n\}$.
\begin{theorem}
\label{uv}
Let $u_\varepsilon(x),v_\varepsilon(x)\in L^2(\Omega)$ two-scale converge to 
$u_0(x,y),v_0(x,y)\in L^2(\Omega\times Y)$ respectively. And also
$\lim\limits_{\varepsilon\to 0}\|u_\varepsilon\|_{L^2(\Omega)}=
\|u_0\|_{L^2(\Omega\times Y)}$
then
$$
\lim\limits_{\varepsilon\to 0}\int_{\Omega}
u_\varepsilon(x)v_\varepsilon(x)\,dx=
\int_{\Omega}\int_{Y}u_0(x,y)v_0(x,y)\,dx\,dy.
$$
\end{theorem}
\begin{proof}
See the proof of Th.1.8 in ~\cite{Allaire}. There one can choose $\psi_n(x,y)$ from 
${\mathcal B}_{TF}={\cal D}(\Omega\times Y)$, $\phi(x)=1$ even if $1$ is not in 
${\cal D}(\Omega)$. Note that $v_\varepsilon$ must be
bounded in assumptions of Th.1.8 and here it is due to (ii) in Def.~\ref{TSC}.   
\qquad
\end{proof}
\begin{corollary}
Necessary conditions (\ref{ATF conditions}) are also sufficient for a function to be ATF.
\end{corollary}
\begin{proof}
In Th.~\ref{uv}, $u_\varepsilon(x)=\phi(x,x/\varepsilon)$. 
$\phi$ satisfies conditions (\ref{ATF conditions}). $v_\varepsilon$ is an arbitrary two-scale convergent
sequence. By Def.~\ref{ATF} $\phi$ is ATF.\qquad   
\end{proof}
With the help of (\ref{ATF conditions}) we can verify whether a particular 
function is ATF. 
The condition (\ref{ATF condition 2}) alone is not enough
\cite[Rem. 1.4.5]{Radu}. Although having a linear space of functions satisfying 
(\ref{ATF condition 2}), there is no need to check (\ref{ATF condition 1}):
\begin{proposition}
Let $L$ be a linear space of functions such that $L\supset{\cal B}_{TF}$ 
and all functions from $L$ satisfy (\ref{ATF condition 2}). 
Then $L\subset{\cal A}_{TF}$.
\end{proposition}
\begin{proof}
We have to check (\ref{ATF condition 1}) for $\phi\in L$. 
$u_\varepsilon(x)=\phi(x,x/\varepsilon)$ is bounded due to 
(\ref{ATF condition 2}). For any $\psi(x,y)\in{\cal B}_{TF}$:
$$
\phi(x,x/\varepsilon)\psi(x,x/\varepsilon)=\frac{1}{2}\left\{
\left[\phi(x,x/\varepsilon)+\psi(x,x/\varepsilon)\right]^2
-\phi(x,x/\varepsilon)^2-\psi(x,x/\varepsilon)^2\right\}
$$
$$
\lim\limits_{\varepsilon\to 0}
\int\limits_{\Omega}\phi(x,\frac{x}{\varepsilon})\psi(x,\frac{x}{\varepsilon})\,dx
=\int\limits_{\Omega}\int\limits_Y\phi(x,y)\psi(x,y)\,dx\,dy
$$
We used (\ref{ATF condition 2}) for $\phi+\psi,\phi,\psi\in L$.\qquad
\end{proof}
\begin{proposition}
${\cal A}_{TF}$ is a linear space.
\end{proposition}
\begin{proof}
Let $\phi_1,\phi_2\in{\cal A}_{TF}$, real numbers
$\alpha,\beta$. We need to check (\ref{ATF conditions}) for 
$\alpha\phi_1+\beta\phi_2$. For any $\psi\in{\cal B}_{TF}$, (\ref{ATF condition 1}) is valid:
$$
\lim\limits_{\varepsilon\to 0}\int\limits_{\Omega}
\left[(\alpha\phi_1+\beta\phi_2)\psi\right](x,x/\varepsilon)\, dx=
\int\limits_{\Omega}\int\limits_Y
\left[(\alpha\phi_1+\beta\phi_2)\psi\right](x,y)\,dx\,dy.
$$
We can use (\ref{ATF condition 2}) for $\phi_1,\phi_2$ and Th.~\ref{uv} with 
$u_\varepsilon(x)=\phi_1(x,x/\varepsilon)$, 
 $v_\varepsilon(x)=\phi_2(x,x/\varepsilon)$ to verify (\ref{ATF condition 2}) for $\alpha\phi_1+\beta\phi_2$:
$$
[\alpha\phi_1+\beta\phi_2]^2(x,x/\varepsilon)=
[\alpha^2\phi_1^2+\beta^2\phi_2^2+2\alpha\beta\phi_1\phi_2](x,x/\varepsilon).
$$
Hence
$$
\lim\limits_{\varepsilon\to 0}\int\limits_{\Omega}
[\alpha\phi_1+\beta\phi_2]^2(x,x/\varepsilon)\, dx=
\int\limits_{\Omega}\int\limits_Y
[\alpha\phi_1+\beta\phi_2]^2(x,y)\,dx\,dy.
$$
\qquad
\end{proof}

The following sections contain properties of the two--scale 
${\mathcal Continuous}$ and ${\mathcal Discrete}$ 
extensions of $a_M$ respectively. 
Our main goal is to show that these extensions $a(x,y)$ are ATF. 
If $a_M\in L^1(\tilde\Omega)$ then we assume that $a(x,y)$ and $a(x,x/\varepsilon)$ are 
constructed pointwise a.e. in $\Omega\times Y$ and in $\Omega$ 
from some representative $a_M(x)$ of $a_M$.
Another representative $\bar a_M(x)$ results in a.e. the same functions
$\bar a(x,y)$ and $\bar a(x,x/\varepsilon)$.
\section{Properties of the two-scale ${\mathcal Continuous}$ Extension}
\label{Section Properties of the two-scale Continuous Extension}  
In this section we deal only with the extension $a(x,y)$ constructed from
$a_M(x)$ in the subsection \ref{Section Continuous Extension}.
\begin{proposition}
\label{P1 B}
For fixed $x\in\Omega$, $a(x,\cdot)$ was constructed piecewise from 
$a_M(\cdot)$, namely $\mathbb{R}^d$ is divided into $1^d$-cubes, by the grid
$$
{\mathcal N}_x(x)=\left\{y\in\mathbb{R}^d\mid \exists k\in\{1,\dots, d\}, 
i\in\mathbb{Z}:\quad y_k=x_k/\bar\varepsilon+i-1/2\right\},
$$
each cube corresponds to the same $\bar\varepsilon$-cube $W(x)$.
\end{proposition}
\begin{proposition}
\label{P2 B}
Let us now fix some $y\in\mathbb{R}^d$. The function $a(\cdot,y)$ is piecewise 
constant on $x\in\Omega$: for each $y$, $\Omega$ is divided by cubic 
$\bar\varepsilon$ grid 
$$
{\mathcal N}_y(y)=\{x\in\Omega\mid y\in {\mathcal N}_x(x)\}=
\left\{x\in\Omega\mid \exists k,i:\quad x_k=y_k\bar\varepsilon-(i-1/2)\bar\varepsilon\right\}
$$
into parts where $a(\cdot,y)$ is constant.
\end{proposition}

The way in which $a(x,y)$ was constructed makes it difficult to deal with
$a(x,x/\varepsilon)$. We need a simple representation of $a(x,x/\varepsilon)$ 
for $\varepsilon\ne\bar\varepsilon$. The first argument $x$ determines the set 
${\mathcal N}_x(x)$ in $\mathbb R^d$. The second argument $x/\varepsilon$ determines
which value of $a_M(x)$ in the neighbourhood $W(x)$ should be taken as the 
value $a(x,x/\varepsilon)$. The non--periodicity of $a_M(x)$ causes an
uncertainty when $x/\varepsilon\in {\mathcal N}_x(x)$.
$$
{\mathcal N}=\left\{x\in\mathbb{R}^d\mid x/\varepsilon\in {\mathcal N}_x(x)\right\}=
\left\{x\in\mathbb{R}^d\mid\exists k,i:\quad
x_k/\varepsilon=x_k/\bar\varepsilon+i-1/2\right\}
$$
$$
{\mathcal N}=\left\{x\in\mathbb{R}^d\mid\exists k\in\{1,\dots,d\},i\in\mathbb{Z}:  
x_k=\left.(i-1/2)\varepsilon\bar\varepsilon\right/
(\bar\varepsilon-\varepsilon)\right\}
$$
${\mathcal N}$ divides $\mathbb{R}^d$ into open cubes $\tilde\Delta_I$ with a side 
$\Delta=\varepsilon\bar\varepsilon/|\bar\varepsilon-\varepsilon|$
and centers in
$$
\dot x_I=\frac{\varepsilon\bar\varepsilon}{\bar\varepsilon-\varepsilon}I,
\qquad I=(i_1\dots i_d)\in \mathbb{Z}^d
$$
Let ${\bf J}_\varepsilon$ be a set of multiindexes $I\in\mathbb{Z}^d$
that $\Delta_I:=\tilde\Delta_I\cap\Omega$ is not an empty set.

If $\dot x_I\in\Omega$ then
$$
a\left(\dot x_I,\frac{\dot x_I}{\varepsilon}\right)=\tilde 
a\left(\dot x_I,\frac{\bar\varepsilon}{\varepsilon}\dot x_I\right)=
\tilde a(\dot x_I,\dot x_I+\bar\varepsilon I)=
\{\mbox{$\bar\varepsilon$-periodicity}\}=
\tilde a(\dot x_I,\dot x_I)=
a_M(\dot x_I) 
$$
Similar if $x\in\Delta_I$ then $x=\dot x_I+h$, $|h_k|<\Delta/2$
$$
a\left(\dot x_I+h,\frac{\dot x_I+h}{\varepsilon}\right)=\tilde
a(\dot x_I+h,\frac{\bar\varepsilon}{\varepsilon}(\dot x_I+h))=
\{\mbox{$\bar\varepsilon$-periodicity}\}=
\tilde a(\dot x_I+h,\dot x_I+\frac{\bar\varepsilon}{\varepsilon}h)=
$$
$$
=a_M(\dot x_I+\frac{\bar\varepsilon}{\varepsilon}h)\qquad
\mbox{since}\qquad
\dot x_I+\frac{\bar\varepsilon}{\varepsilon}h\in W(\dot x_I+h).
$$
We have proved the following 
\begin{proposition}
\label{P3 B}
The simple representation of $a(x,x/\varepsilon)$ for all $\varepsilon>0$, 
$\varepsilon\ne\bar\varepsilon$ is:
$$
\mbox{if $x\in\Delta_I$\qquad then}\qquad
a\left(x,\frac{x}{\varepsilon}\right)=
a_M(\dot x_I+\frac{\bar\varepsilon}{\varepsilon}(x-\dot x_I))
$$
or using the Heaviside function ${\bf 1}_{\Delta_I}(x)$ being $1$ in 
$\Delta_I$ and $0$ elsewhere we have: 
\begin{equation}
\label{a(x,x/epsilon) Heaviside}
a(x,x/\varepsilon)=\sum\limits_{I\in{\bf J}_\varepsilon}
{\bf 1}_{\Delta_I}(x)a_M\bigl(\dot x_I+\frac{\bar\varepsilon}
{\varepsilon}(x-\dot x_I)\bigr)
\end{equation}
\end{proposition}
Roughly speaking for $\varepsilon<\bar \varepsilon$ 
[$\varepsilon>\bar \varepsilon$] 
$a(x,x/\varepsilon)$ is built from compressed [stretched] 
cubes taken from $a_M(x)$.

{\it For the following let $\phi(x,y)$ be a function from 
$C(\overline\Omega\times\overline Y)$, $Y$-periodic in $y$. 
We will consider $a(x,y)\phi(x,y)$. 
Important particular  case: $\phi(x,y)=1$.}

\begin{proposition}
\label{P4 B} 
1) If $a_M(x)$ is measurable in $\tilde\Omega$ then 
$a(x,x/\varepsilon)\phi(x,x/\varepsilon)$ 
is measurable in $\Omega$. 2) If $a_M(x)\in L^1(\tilde\Omega)$, then
$a(x,x/\varepsilon)\phi(x,x/\varepsilon)\in L^1(\Omega)$.
\end{proposition}
\begin{proof}
We only have to consider the case $\varepsilon\ne\bar\varepsilon$.\\
1) $\phi(x,x/\varepsilon)\in C(\overline\Omega)$ is measurable.
$a(x,x/\varepsilon)$ is measurable since it is a sum of measurable functions
(\ref{a(x,x/epsilon) Heaviside}).\\ 
2) 
If $\Omega$ is bounded with the diameter $2R$, then 
$\bigcup_{I\in{\bf J}_\varepsilon}\tilde\Delta_I$ is bounded with the diameter 
$D:=2(R+\sqrt{d}\Delta)$, 
$$
\sum_{I\in{\bf J}_\varepsilon}\Delta^d=
\sum_{I\in{\bf J}_\varepsilon}\mu(\tilde\Delta_I)\le D^d\quad 
\Rightarrow\quad
\sum_{I\in{\bf J}_\varepsilon}\frac{\varepsilon^d\bar\varepsilon^d}
{|\bar\varepsilon-\varepsilon|^d}\le D^d\quad
\Rightarrow\quad
\sum_{I\in{\bf J}_\varepsilon}\frac{\varepsilon^d}{\bar\varepsilon^d}
\le\frac{|\bar\varepsilon-\varepsilon|^d}{\bar\varepsilon^{2d}} D^d, 
$$
$$
\widetilde W(\dot x_I):=
\left\{z=\dot x_I+\frac{\bar\varepsilon}{\varepsilon}(x-\dot x_I)\mid
x\in\Delta_I\right\}=
\left\{z=\dot x_I+\frac{\bar\varepsilon}{\varepsilon}h\mid 
\dot x_I+h\in\Delta_I\right\}\subset{\tilde\Omega}
$$ 
$$
\|a(x,x/\varepsilon)\|_{L^1(\Omega)}=
\int\limits_\Omega |a(x,x/\varepsilon)|\, dx=
\sum\limits_{I\in{\bf J}_\varepsilon}
\int\limits_{\Delta_I}
|a_M(\dot x_I+\frac{\bar\varepsilon}{\varepsilon}(x-\dot x_I))|\, dx\le
$$
$$
\le\sum\limits_{I\in{\bf J}_\varepsilon}
\frac{\varepsilon^d}{\bar\varepsilon^d}
\int\limits_{\widetilde W(\dot x_I)}|a_M(z)|\, dz
\le\|a_M\|_{L^1(\tilde\Omega)}
\sum_{I\in{\bf J}_\varepsilon}\frac{\varepsilon^d}{\bar\varepsilon^d}
\le\|a_M\|_{L^1(\tilde\Omega)}
\frac{|\bar\varepsilon-\varepsilon|^d}{\bar\varepsilon^{2d}}D^d. 
$$
$$
\int\limits_\Omega |a(x,x/\varepsilon)\phi(x,x/\varepsilon)|\,dx\le
\|\phi\|_C\int\limits_\Omega |a(x,x/\varepsilon)|\,dx\le
\|\phi\|_C\|a_M\|_{L^1(\tilde\Omega)}
\frac{|\bar\varepsilon-\varepsilon|^d}{\bar\varepsilon^{2d}}D^d.
$$
\end{proof}
\begin{proposition}
\label{P6 B}
If $a_M(x)\in L^1(\tilde\Omega)$ then
$$
M(x)=\int_Y a(x,y)\phi(x,y)\, dy,\quad
M_+(x)=\int_Y |a(x,y)\phi(x,y)|\, dy 
\quad\mbox{are continuous in $\overline\Omega$,}
$$
$M(x)$,$M_+(x)$ are bounded by
$\|\phi\|_C\|a_M\|_{L^1(\tilde\Omega)}/\bar\varepsilon^d$.
\end{proposition}
\begin{proof}
$a(x,\cdot)\in L^1(Y)$ since it was constructed from $a_M(\cdot)$. Therefore
$M(x)$ and $M_+(x)$ are well defined.
To show continuity let us fix an arbitrary $E>0$.
$$
|M(x+h)-M(x)|=\Bigl|\int_Y a(x+h,y)\phi(x+h,y)\, dy
-\int_Ya(x,y)\phi(x,y)\, dy\Bigr|\le
$$
$$
\Bigl|\int\limits_Y a(x+h,y)[\phi(x+h,y)-\phi(x,y)]\, dy\Bigr|+
\Bigl|\int\limits_Y a(x+h,y)\phi(x,y)\,dy-\int\limits_Y a(x,y)\phi(x,y)\,dy\Bigr|
$$
For continuous $\phi$ one can find such $\delta_1$ that 
$|\phi(x+h,y)-\phi(x,y)|<E\bar\varepsilon^d/2\|a_M\|_{L^1(\tilde\Omega)}$
when $|h|_\infty<\delta_1$ ($|h|_\infty=\max\limits_k |h_k|$ we distinguish from the
vector's 
absolute value $|h|=\sqrt{\sum_k h_k^2}$). 
This means that the first absolute value is less than $E/2$.\\ 
Now we consider the second absolute value. 
Using that $a(x,y)=\tilde a(x,\bar\varepsilon y)$, 
substitution of variables $z=\bar\varepsilon y$ we obtain
$$
\int_Y a(x,y)\phi(x,y)\,dy=\frac{1}{\bar\varepsilon^d}
\int_{\bar\varepsilon Y}\tilde a(x,z)\phi(x,z/\bar\varepsilon)\,dz=
$$
$\tilde a(x,z)\phi(x,z/\bar\varepsilon)$ is $\bar\varepsilon$-periodic in $z$,
integral over $\bar\varepsilon Y$ is  equal to integral over any 
$\bar\varepsilon$ cube
$$
=\frac{1}{\bar\varepsilon^d}
\int_{W(x)}\tilde a(x,z)\phi(x,z/\bar\varepsilon)\,dz
=\frac{1}{\bar\varepsilon^d}\int_{W(x)}a_M(z)\phi(x,z/\bar\varepsilon)\,dz
$$
Similar
$$
\int_Y a(x+h,y)\phi(x,y)\,dy=
\frac{1}{\bar\varepsilon^d}\int_{W(x+h)}a_M(z)\phi(x,z/\bar\varepsilon)\,dz
$$
$$
\frac{1}{\bar\varepsilon^d}\Bigl|
\int\limits_{W(x+h)}a_M(z)\phi(x,z/\bar\varepsilon)\,dz-
\int\limits_{W(x)}a_M(z)\phi(x,z/\bar\varepsilon)\,dz\Bigr|\le
\frac{\|\phi\|_C}{\bar\varepsilon^d}
\int\limits_{W(x+h)\triangle W(x)}|a_M(z)|\,dz 
$$
$\mu\left(W(x+h)\triangle W(x)\right)\le 2d\bar\varepsilon^{d-1}|h|_{\infty}$
Using absolute continuity of Lebesgue integral, there 
exists $\delta_2$: $|h|_{\infty}<\delta_2$ guarantees that the second absolute
value is less than $E/2$ and consequently for 
$|h|_{\infty}<\min\{\delta_1,\delta_2\}$ we have  $|M(x+h)-M(x)|<E$.
$$
2)\qquad |M(x)|\le\frac{\|\phi\|_C}{\bar\varepsilon^d}
\int_{\tilde\Omega}|a_M(x)|\, dx=
\left.\|\phi\|_C\|a_M\|_{L^1(\tilde\Omega)}\right/\bar\varepsilon^d.
$$
Similar with $M_+(x)$.
\end{proof}
\begin{proposition}
\label{P7 B} 
If $a_M(x)$ is measurable in $\tilde\Omega$ then $a(x,y)\phi(x,y)$ 
is measurable in $\Omega\times Y$.
\end{proposition}
\begin{proof}
$\phi(x,y)$ is continuous hence measurable.
To show measurability of $a(x,y)$ we will construct a sequence of 
measurable functions $\{a_\delta(x,y)\}$ converging to $a(x,y)$ a.e when
$\delta\to 0$.
Let us divide $\mathbb{R}^d$ into cubes 
$\Box^\delta_i=\Bigl[i_1\delta,(i_1+1)\delta\Bigr)\times
\dots\times\Bigl[i_d\delta,(i_d+1)\delta\Bigr)$, 
$i\in\mathbb{Z}^d$. ${\bf I}_\delta$ is a set of indexes $i\in\mathbb{Z}^d$ that 
$\Box^\delta_i\cap\Omega\ne\emptyset$. $\delta$ is small enough that 
$\Omega\subset\bigcup\limits_{i\in I_\delta}\Box^\delta_i\subset\tilde\Omega$.
For $i\in{\bf I}_\delta$ let $\tilde x^\delta_i$ be an arbitrary point of 
$\Box^\delta_i\cap\Omega$ (e.g. the center).
$$
\mbox{The function}\qquad
a_\delta(x,y):=a(\tilde x^\delta_i,y),\qquad
\mbox{when $x\in\Box^\delta_i\cap\Omega$, $i\in I_\delta$},
$$
is measurable in $\Omega\times Y$ since $a(\tilde x^\delta_i,y)$ is a 
measurable function in $Y=(0,1)^d$ and $\Box^\delta_i\cap\Omega$ is a measurable 
set. We have to show that the sequence pointwise converges to $a(x,y)$ in 
$\Omega\times Y\setminus O$, where 
$O:=\{(x,y)\in\Omega\times Y\mid y\in {\mathcal N}_x(x)\}$ is a zero measure set.
$$
O\cap\left(\Box^\delta_i\times Y\right)\subset O^\delta_i:=
\Box^\delta_i\times
\Bigl(Y\cap\bigcup\limits_{x\in\Box^\delta_i}{\mathcal N}_x(x)\Bigr),\qquad
O\subset\bigcup\limits_{i\in{\bf I}_\delta}O^\delta_i.
$$
$O^\delta_i$ is a measurable set,
$\mu_{X\times Y}(O^\delta_i)=\mu_X(\Box^\delta_i)\times
\mu_Y(Y\cap\bigcup\limits_{x\in\Box^\delta_i}{\mathcal N}_x(x))
\le \mu_X(\Box^\delta_i)d
\frac{\delta}{\bar\varepsilon}
$. 
$$
\mu_{X\times Y}(O)\le\sum\limits_{i\in{\bf I}_\delta}d\frac{\delta}{\bar\varepsilon}
\mu_X(\Box^\delta_i)\le d\frac{\delta}{\bar\varepsilon}\mu_X(\tilde\Omega)\to 0,
\mbox{ when $\delta\to 0$}\quad\Rightarrow\quad\mu_{X\times Y}(O)=0.
$$
Let $(x,y)\in(\Omega\times Y)\setminus O$. 
It means that $dist(y,{\mathcal N}_x(x))>0$ (here $dist(y,\hat y)=|y-\hat y|_\infty$). 
If we consider a $\delta$-partition with 
$\delta<\bar\varepsilon\, dist(y,{\mathcal N}_x(x))$, 
$x\in\Box^\delta_i$ for some $i$ then $(x,y)\notin O^\delta_i$ since
for all $\hat x\in\Box^\delta_i$, $y$ is enough far from ${\mathcal N}_x(\hat x)$.
As we know from Prop.~\ref{P2 B} $a(\cdot,y)$ is piecewise constant
in $\Omega$ and it changes value at those $\hat x$ that $y\in {\mathcal N}_x(\hat x)$.
The whole set $\Box^\delta_i\cap\Omega$ belongs to the cube where $a(\cdot,y)$ is constant.
As a result: $\forall\hat x\in\Box^\delta_i\cap\Omega$, 
$a(\hat x,y)=a(\tilde x^\delta_i,y)$.  
On the other hand from the definition 
of $a_\delta$: $\forall\hat x\in\Box^\delta_i\cap\Omega$, 
$a_\delta(\hat x,y)=a(\tilde x^\delta_i,y)$. 
Consequently for our particular point $(x,y)\in(\Box^\delta_i\times Y)\setminus O$
and small enough $\delta$ we have $a_\delta(x,y)=a(x,y)$.
\end{proof}
\begin{lemma}
\label{P8 B}
Let $a_M(x)\in L^1(\tilde\Omega)$; $\phi(x,y)\in 
C(\overline\Omega\times\overline Y)$, $Y$-periodic in $y$; 
$a(x,y)$ is the ${\mathcal Continuous}$ Extension of $a_M(x)$.
Then
\begin{equation}
\lim\limits_{\varepsilon\to 0}
\int_\Omega a(x,x/\varepsilon)\phi(x,x/\varepsilon)\, dx=
\int_\Omega\int_Y a(x,y)\phi(x,y)\,dx\,dy.
\label{admissible}
\end{equation}
\end{lemma}
\begin{proof}
1. $a(x,y)\phi(x,y)$ is measurable (Prop.~\ref{P7 B}). 
$M_+(x)$ is continuous 
(Prop.~\ref{P6 B}).
From Fubini's theorem $a(x,y)\phi(x,y)\in L^1(\Omega\times Y)$, 
the right hand side is well defined.\\
2. Integrals in the left hand side are well defined (Prop.~\ref{P4 B}).\\
3. Now we need to check the equality (\ref{admissible}).
Let $E>0$ be an arbitrarily small number,
for some $\varepsilon$ ($\varepsilon\le\bar\varepsilon/2$) we consider a
subdivision of $\Omega$ (already defined in Prop.~\ref{P3 B}) with 'central' 
points $\dot x_I$:
$$
\overline{\Omega}=\bigcup_{I\in{\bf J}_\varepsilon}\overline{\Delta_I};\quad  
{\bf J}_\varepsilon={\bf J}_\varepsilon^{int}\cup({\bf J}_\varepsilon\setminus
{\bf J}_\varepsilon^{int});\quad
{\bf J}_\varepsilon^{int}=
\{I\in{\bf J}_{\varepsilon}\mid\tilde\Delta_I\subset\Omega\};
\quad \overline{\Omega^{int}}=
\bigcup_{I\in{\bf J}_\varepsilon^{int}}\overline{\Delta_I}. 
$$
For not too bad $\partial\Omega$ and small enough $\varepsilon$, 
$\mu(\Omega\setminus\Omega^{int})$ is arbitrarily small:
$$
\mu(\Omega\setminus\Omega^{int})\le\sum_{I\in{\bf J}_\varepsilon\setminus
{\bf J}_\varepsilon^{int}}\Delta^d
\le\frac{\bar\varepsilon^d E}{5\|\phi\|_C\|a_M\|_{L^1(\tilde\Omega)}}.
$$
We will approximate the integrals over $\Omega$ using the integrals over 
$\Omega^{int}$. Let us estimate the errors in a similar way as it was done in 
Prop.~\ref{P4 B}, Prop.~\ref{P6 B}:
$$
\int_{\Omega\setminus\Omega^{int}}|a(x,x/\varepsilon)\phi(x,x/\varepsilon)|\,dx
\le\|\phi\|_C\sum_{I\in{\bf J}_\varepsilon\setminus
{\bf J}_\varepsilon^{int}}\frac{\varepsilon^d}{\bar\varepsilon^d}
\int_{\widetilde W(\dot x_I)}|a_M(z)|\, dz\le
$$
$$
\le\|\phi\|_C\|a_M\|_{L^1(\tilde\Omega)}\frac{|\bar\varepsilon-\varepsilon|^d}
{\bar\varepsilon^{2d}}
\frac{\bar\varepsilon^d E}{5\|\phi\|_C\|a_M\|_{L^1(\tilde\Omega)}}
\le E/5,\qquad
\int_{\Omega\setminus\Omega^{int}}|M(x)|\, dx\le E/5.
$$
Since $M(x)$ is continuous in $\overline\Omega$, 
then for small enough $\varepsilon$, the right hand side integral in 
(\ref{admissible}) can be approximated by the sum
$$
\int_{\Omega}\int_Y a(x,y)\phi(x,y)\, dy\,dx\approx
\sum_{I\in{\bf J}_\varepsilon^{int}}\mu(\Delta_I)
\int_Y a(\dot x_I,y)\phi(\dot x_I,y)\, dy=
$$
\begin{equation}
=\sum_{I\in{\bf J}_\varepsilon^{int}}\frac{\mu(\Delta_I)}{\bar\varepsilon^d}
\int_{\bar\varepsilon Y}\tilde a(\dot x_I,z)
\phi(\dot x_I,\frac{z}{\bar\varepsilon})\, dz=
\sum_{I\in{\bf J}_\varepsilon^{int}}\frac{\mu(\Delta_I)}{\mu(W(\dot x_I))}
\int_{W(\dot x_I)}a_M(z)\phi(\dot x_I,\frac{z}{\bar\varepsilon})\, dz.
\label{lhs appr}
\end{equation}
with error not greater than $2E/5$.
  
The integral in the left hand side of (\ref{admissible})
can be approximated by
$$
\int_\Omega a(x,x/\varepsilon)\phi(x,x/\varepsilon)\, dx\approx
\sum_{I\in{\bf J}_\varepsilon^{int}}\int_{\Delta_I}
a_M\Bigl(\dot x_I+\frac{\bar\varepsilon}{\varepsilon}(x-\dot x_I)\Bigr)
\phi(x,\frac{x}{\varepsilon}), dx=
$$
or in the new variables $z=\dot x_I+(x-\dot x_I)\bar\varepsilon/\varepsilon$,\quad
$x=\dot x_I+(z-\dot x_I)\varepsilon/\bar\varepsilon$,\quad 
$\Delta_I\to \widetilde W(\dot x_I)$:
$$
=\sum_{I\in{\bf J}_\varepsilon^{int}}
\frac{\varepsilon^d}{\bar\varepsilon^d}
\int_{\widetilde W(\dot x_I)}a_M(z)\phi\Bigl(x(z),
\frac{\dot x_I}{\varepsilon}+\frac{z-\dot x_I}{\bar\varepsilon}\Bigr)\, dz=
$$
\begin{equation}
\label{Initial approximation}
=\sum_{I\in{\bf J}_\varepsilon^{int}}
\frac{\mu(\Delta_I)}{\mu(\widetilde W(\dot x_I))}
\int_{\widetilde W(\dot x_I)}a_M(z)
\phi\Bigl(x(z),\frac{z}{\bar\varepsilon}\Bigr)\, dz.
\end{equation}
The last equality is since $\dot x_I(1/\varepsilon-1/\bar\varepsilon)=I$,
$\phi(x,y)$ is $Y$-periodic and 
$\varepsilon/\bar\varepsilon=\mu(\Delta_I)/\mu(\widetilde W(\dot x_I))$. 
The approximation error is not greater than $E/5$. 
Further approximation of (\ref{Initial approximation}):
\begin{equation}
\approx\sum_{I\in{\bf J}_\varepsilon^{int}}
\frac{\mu(\Delta_I)}{\mu(\widetilde W(\dot x_I))}
\int_{\widetilde W(\dot x_I)}a_M(z)
\phi\Bigl(\dot x_I,\frac{z}{\bar\varepsilon}\Bigr)\, dz
\label{rhs appr}
\end{equation}
has an error 
$$
\sum_{I\in{\bf J}_\varepsilon^{int}}
\frac{\mu(\Delta_I)}{\mu(\widetilde W(\dot x_I))}
\int_{\widetilde W(\dot x_I)}a_M(z)\left[
\phi\Bigl(x(z),\frac{z}{\bar\varepsilon}\Bigr)-
\phi\Bigl(\dot x_I,\frac{z}{\bar\varepsilon}\Bigr)\right]\, dz
$$
which we can estimate in absolute value like in Prop.~\ref{P4 B}
(restricting to $\varepsilon\le\bar\varepsilon/2$):
$$
\delta\|a_M\|_{L^1(\tilde\Omega)}\sum_{I\in{\bf J}_\varepsilon^{int}}
\frac{\varepsilon^d}{\bar\varepsilon^d}\le
\delta\|a_M\|_{L^1(\tilde\Omega)}\frac{|\bar\varepsilon-\varepsilon|^d}
{\bar\varepsilon^{2d}}\mu(\Omega)\le 
\delta\frac{\|a_M\|_{L^1(\tilde\Omega)}}{\bar\varepsilon^d}\mu(\Omega)\le E/5, 
$$
where $|\phi(x_1,y)-\phi(x_2,y)|<\delta=
E\bar\varepsilon^d/5\|a_M\|_{L^1(\tilde\Omega)}\mu(\Omega)$ when 
$|x_1-x_2|_\infty\le\Delta/2\le\varepsilon$.

Now we compare (\ref{lhs appr}) and (\ref{rhs appr}).
$$
\mbox{For $I\in{\bf J}_\varepsilon^{int}$:}\qquad
\widetilde W(\dot x_I)=\{x\in\mathbb{R}^d\mid |x-\dot x_I|_{\infty}\le
\frac{\bar\varepsilon^2}{2|\bar\varepsilon-\varepsilon|}\},\qquad
W(\dot x_I)\subset\widetilde W(\dot x_I).
$$
$\mu(\widetilde W\setminus W)=\mu(\widetilde W)-\mu(W)=
\bar\varepsilon^{2d}/(\bar\varepsilon-\varepsilon)^d-\bar\varepsilon^d
\stackrel{\varepsilon\to 0}{\longrightarrow}0$. Then for small enough 
$\varepsilon$
$$
\left|\frac{1}{\mu(\widetilde W)}\int_{\widetilde W}
a_M(z)\phi\Bigl(\dot x_I,\frac{z}{\bar\varepsilon}\Bigr)\, dz-
\frac{1}{\mu(W)}\int_{W}
a_M(z)\phi\Bigl(\dot x_I,\frac{z}{\bar\varepsilon}\Bigr)\, dz\right|\le
$$
$$
\le\frac{\mu(\widetilde W)-\mu(W)}{\mu(W)\mu(\widetilde W)}\int_W 
|a_M(z)\phi\Bigl(\dot x_I,\frac{z}{\bar\varepsilon}\Bigr)|\, dz
+\frac{1}{\mu(\widetilde W)}\int_{\widetilde W\setminus W}
|a_M(z)\phi\Bigl(\dot x_I,\frac{z}{\bar\varepsilon}\Bigr)|\, dz\le
$$
$$
\le(\mu(\widetilde W)-\mu(W))
\frac{\|\phi\|_C\|a_M\|_{L^1(\tilde\Omega)}}
{\bar\varepsilon^{2d}}
+\frac{\|\phi\|_C}{\mu(\widetilde W)}
\int_{\widetilde W\setminus W}
|a_M(z)|\, dz\le E/5\mu(\Omega). 
$$
For the second term we used the absolute continuity of Lebesgue integral.

For small enough $\varepsilon$:
$$
|\int_\Omega a(x,x/\varepsilon)\phi(x,x/\varepsilon)\, dx-
\int_\Omega\int_Y a(x,y)\phi(x,y)\,dy\,dx|\le
\frac{4E}{5}+\frac{E}{5\mu(\Omega)}
\sum_{I\in{\bf J}_\varepsilon^{int}}\mu(\Delta_I)\le E.
$$
\end{proof}
\section{Properties of the two-scale ${\mathcal Discrete}$ Extension}
\label{Section Properties of the two-scale Discrete Extension}
In this section we deal only with the extension $a(x,y)$ constructed from
$a_M(x)$ in the Subsection~\ref{Section Discrete Extension}.
\begin{proposition}
\label{P1 C} 
For fixed $x\in\Omega\cap\Omega_j$, $a(x,\cdot)$ was constructed piecewise from 
$a_M(\cdot)$, namely $\mathbb{R}^d$ is divided into $1^d$ cubes by the grid
$$
{\mathcal N}_x(x)=\left\{y\in\mathbb{R}^d\mid\exists k\in\{1\dots d\}, 
i\in\mathbb{Z}:\quad y_k=\hat x^j_k/\bar\varepsilon+i-1/2\right\},
$$
each cube corresponds to the same $\bar\varepsilon$-cube $W_j$.
If $x_1,x_2\in\Omega_j$ then ${\mathcal N}_x(x_1)={\mathcal N}_x(x_2)$.
\end{proposition}
\begin{proposition}
\label{P2 C}
Let us now fix some $y\in\mathbb{R}^d$. The function $a(\cdot,y)$ is a 
constant in each $\Omega_j$.
\end{proposition}

Like in the previous section, we need a simple representation of $a(x,x/\varepsilon)$.
For each $W_j$ if is convenient to correspond an $\bar\varepsilon$-periodic 
function 
$$
\tilde a_j(y):=a_M(y)\qquad\mbox{when $y\in W_j$. Periodically expanded to 
$\mathbb{R}^d$.}
$$
$a_j(y):=\tilde a_j(\bar\varepsilon y)$ is $1^d$-periodic in $\mathbb{R}^d$.

The simple representation of $a(x,y)$ is
\begin{equation}
\label{simple form A(x,y) C}
a(x,y)=\sum_{j=1}^N{\bf 1}_{\Omega_j}(x)\, a_j(y),\qquad\mbox{then}\qquad
a(x,x/\varepsilon)=\sum_{j=1}^N{\bf 1}_{\Omega_j}(x)\, a_j(x/\varepsilon).
\end{equation}
To describe behaviour of $a_j(x/\varepsilon)$
inside $\Omega_j$ let us define the following $\varepsilon$-cubes in 
$\mathbb{R}^d$
$$
\mbox{with centers in 
$\dot x^j_I:=\hat x^j\varepsilon/\bar\varepsilon+I\varepsilon$}
\qquad\tilde\Delta_I^j=\{x\in\mathbb{R}^d\mid |x-\dot x^j_I|_\infty
<\varepsilon/2\}
$$
If $x\in\tilde\Delta_I^j$ then $x=\dot x^j_I+h$
($|h|_\infty<\varepsilon/2$),
$$
a_j(x/\varepsilon)=\tilde a_j(\hat x^j+\bar\varepsilon
I+h\bar\varepsilon/\varepsilon)=\{\mbox{$\bar\varepsilon$-periodicity}\}=
\tilde a_j(\hat x^j+h\bar\varepsilon/\varepsilon)=
a_M(\hat x^j+h\bar\varepsilon/\varepsilon)
$$ 
since $\hat x^j+h\bar\varepsilon/\varepsilon\in W_j$.
We also have a partition of $\Omega_j$:
$$
\overline{\Omega}_j=\bigcup_{I\in{\bf J}_\varepsilon(j)}\overline{\Delta}_I^j,\qquad
\mbox{$I$ belongs to ${\bf J}_\varepsilon(j)$ when 
$\Delta_I^j:=\tilde\Delta_I^j\cap\Omega_j\ne\emptyset$}.
$$
We have proved the following
\begin{proposition}
\label{P3 C}
The simple representation of $a(x,x/\varepsilon)$ for all $\varepsilon>0$ is :
$$
a(x,x/\varepsilon)=\sum_{j=1}^N\sum_{I\in{\bf J}_\varepsilon(j)}
{\bf 1}_{\Delta_I^j}(x)\, 
a_M(\hat x^j+\frac{\bar\varepsilon}{\varepsilon}(x-\dot x^j_I))
$$
\end{proposition}
We again assume that $\phi(x,y)\in C(\overline\Omega\times\overline Y)$, $Y$ --periodic in
$y$.
\begin{proposition}
\label{P4 C}
1)If $a_M(x)$ is measurable in $\tilde\Omega$ then 
$a(x,x/\varepsilon)\phi(x,x/\varepsilon)$
is measurable in $\Omega$. 2) If $a_M(x)\in L^1(\tilde\Omega)$ then
$a(x,x/\varepsilon)\phi(x,x/\varepsilon)\in L^1(\Omega)$.\\
\end{proposition}
\begin{proof}
1) $a_j(x/\varepsilon)$ are measurable; $\Omega_j$ are measurable sets; 
$a(x,x/\varepsilon)$ is a sum of measurable functions. 
$\phi(x,x/\varepsilon)\in C(\overline\Omega)$ is measurable.
$$
\mbox{2) }\int_\Omega |a(x,x/\varepsilon)|\le
\sum_{j=1}^N\sum_{I\in{\bf J}_\varepsilon(j)}\int_{\tilde\Delta_I^j}
|a_M(\hat x^j+\frac{\bar\varepsilon}{\varepsilon}(x-\dot x^j_I))|=
\sum_{j=1}^N\sum_{I\in{\bf J}_\varepsilon(j)}
\left(\frac{\varepsilon}{\bar\varepsilon}\right)^d
\int_{W_j}|a_M|\le
$$
$$
\le\frac{\|a_M\|_{L^1(\tilde\Omega)}}{\bar\varepsilon^d}\sum_{j=1}^N
\sum_{I\in{\bf J}_\varepsilon(j)}\varepsilon^d
=\frac{\|a_M\|_{L^1(\tilde\Omega)}}{\bar\varepsilon^d}\sum_{j=1}^N
\mu\Bigl(\bigcup_{I\in{\bf J}_\varepsilon(j)}\tilde\Delta^j_I\Bigr)
\le N\|a_M\|_{L^1(\tilde\Omega)}
\left(\frac{\bar\varepsilon+2\varepsilon}{\bar\varepsilon}\right)^d
$$
the measures were estimated by $(\bar\varepsilon+2\varepsilon)^d$ since 
$\Omega_j\subset W_j$ and $W_j$, $\tilde\Delta^j_I$ have sides 
$\bar\varepsilon$, $\varepsilon$ respectively. Therefore
$$
\int_\Omega |a(x,x/\varepsilon)\phi(x,x/\varepsilon)|\,dx\le
\|\phi\|_C\int_\Omega |a(x,x/\varepsilon)|\,dx\le
N\|\phi\|_C\|a_M\|_{L^1(\tilde\Omega)}
\left(\frac{\bar\varepsilon+2\varepsilon}{\bar\varepsilon}\right)^d.
$$
\end{proof}
\begin{proposition}
\label{P6 C}
If $a_M(x)\in L^1(\tilde\Omega)$ then
$$
M(x)=\int_Y a(x,y)\phi(x,y)\, dy,\quad
M_+(x)=\int_Y |a(x,y)\phi(x,y)|\, dy 
$$
are continuous in each $\Omega_j$
and $M(x)$,$M_+(x)$ are bounded by
$\left.\|\phi\|_C\|a_M\|_{L^1(\tilde\Omega)}\right/\bar\varepsilon^d$.
\end{proposition}
\begin{proof}
$$
M(x)=\int_Y a(x,y)\phi(x,y)\, dy=\int_Y a_j(y)\phi(x,y)\, dy=
\frac{1}{\bar\varepsilon^d}
\int_{\bar\varepsilon Y}\tilde a_j(z)\phi(x,z/\bar\varepsilon)\, dz=
$$
$$
=\{\mbox{$\tilde a_j(z)\phi(x,z/\bar\varepsilon)$  is 
$\bar\varepsilon$-periodic in $z$}\}=
\frac{1}{\mu(W_j)}\int_{W_j}a_M(z)\phi(x,z/\bar\varepsilon)\, dz.
$$
Let $E>0$ be an arbitrarily small number.
For $x_1,x_2\in\Omega_j$, $|\phi(x_1,y)-\phi(x_2,y)|\le\delta=
E\bar\varepsilon^d/\|a_M\|_{L^1(\tilde\Omega)}$, 
when
$|x_1-x_2|_{\infty}$ is small enough:
$$
|M(x_1)-M(x_2)|\le\frac{1}{\mu(W_j)}\int_{W_j}|a_M(z)|
|\phi(x_1,z/\bar\varepsilon)-\phi(x_2,z/\bar\varepsilon)|\, dz\le
E. 
$$
$$
|M(x)|\le\frac{\|\phi\|_C}{\bar\varepsilon^d}\int_{W_j}|a_M(z)|\,dz
\le\frac{\|\phi\|_C}{\bar\varepsilon^d}\int_{\tilde\Omega}|a_M(z)|\,dz. 
$$
Similar with $M_+(x)$.
\end{proof}
\begin{proposition}
\label{P7 C}
If $a_M(x)$ is measurable in $\tilde\Omega$ then $a(x,y)\phi(x,y)$ 
is measurable in $\Omega\times Y$.
\end{proposition}
\begin{proof}
$a(x,y)$ is a sum of measurable 
functions (\ref{simple form A(x,y) C}), $\phi(x,y)$ is measurable. 
\end{proof}
\begin{lemma}
\label{P8 C} 
Let $a_M(x)\in L^1(\tilde\Omega)$; 
$\phi(x,y)\in C(\overline\Omega\times\overline Y)$, $Y$-periodic in $y$; 
$a(x,y)$ is the ${\mathcal Discrete}$ Extension of $a_M(x)$.
Then
\begin{equation}
\lim\limits_{\varepsilon\to 0}\int_\Omega 
a(x,x/\varepsilon)\phi(x,x/\varepsilon)\, dx=
\int_\Omega\int_Y a(x,y)\phi(x,y)\,dx\,dy.
\label{admissible 2}
\end{equation}
\end{lemma}
\begin{proof}
Left and right integrals are well defined. Let us fix a small $E>0$.
\begin{equation}
\label{P8C RHS}
\int_\Omega\int_Y a(x,y)\phi(x,y)\,dx\,dy=\sum_{j=1}^N\frac{1}{\mu(W_j)}
\int_{W_j}a_M(z)\Bigl[\int_{\Omega_j}\phi(x,z/\bar\varepsilon)\, dx\Bigr]\, dz.
\end{equation}
Let ${\bf J}_\varepsilon^{int}(j)$ consists of those 
$I\in{\bf J}_\varepsilon(j)$ that $\tilde\Delta^j_I\subset\Omega_j$.
$\overline{\Omega}_j^{int}=\bigcup_{I\in{\bf J}_\varepsilon^{int}(j)}\overline{\Delta}^j_I$.
For not too bad boundaries $\partial\Omega_j$ one can find such small $\varepsilon$ that
$$
\sum\limits_{j=1}^N\mu(\Omega_j\setminus\Omega_j^{int})=
\sum\limits_{j=1}^N\sum_{I\in{\bf J}_\varepsilon(j)\setminus
{\bf J}_\varepsilon^{int}(j)}\varepsilon^d\le
\frac{E\bar\varepsilon^d}{2\|\phi\|_C\|a_M\|_{L^1(\tilde\Omega)}}.
$$
Therefore the error of the approximation
$$
\int_{\Omega}a(x,x/\varepsilon)\phi(x,x/\varepsilon)\, dx=
\sum_{j=1}^N\sum_{I\in{\bf J}_\varepsilon(j)}\int_{\Delta_I^j}
a_M(\hat x^j+\frac{\bar\varepsilon}{\varepsilon}(x-\dot x^j_I))
\phi(x,x/\varepsilon)\, dx\approx
$$
\begin{equation}
\label{P8C approximation 1}
\approx\sum_{j=1}^N\sum_{I\in{\bf J}_\varepsilon^{int}(j)}\int_{\Delta_I^j}
a_M(\hat x^j+\frac{\bar\varepsilon}{\varepsilon}(x-\dot x^j_I))
\phi(x,x/\varepsilon)\, dx
\end{equation}
can be estimated in absolute value as
$$
\sum_{j=1}^N
\sum_{I\in{\bf J}_\varepsilon(j)\setminus{\bf J}_\varepsilon^{int}(j)}
\|\phi\|_C\int_{\Delta_I^j}
|a_M(\hat x^j+\frac{\bar\varepsilon}{\varepsilon}(x-\dot x^j_I))|\, dx\le
$$
$$
\le\sum_{j=1}^N
\sum_{I\in{\bf J}_\varepsilon(j)\setminus{\bf J}_\varepsilon^{int}(j)}
\frac{\|\phi\|_C}{\bar\varepsilon^d}\varepsilon^d\int\limits_{W_j}|a_M(z)|\, dz
\le\frac{\|\phi\|_C\|a_M\|_{L^1(\tilde\Omega)}}{\bar\varepsilon^d}
\sum_{j=1}^N\sum_{I\in{\bf J}_\varepsilon(j)\setminus{\bf
J}_\varepsilon^{int}(j)}\varepsilon^d\le\frac{E}{2}.
$$
We continue (\ref{P8C approximation 1}) by introducing  
new variables for each integral over $\Delta_I^j$: 
$z=\hat x^j+\frac{\bar\varepsilon}{\varepsilon}(x-\dot x^j_I)$,
$x=x^j_I(z)=\dot x^j_I+\frac{\varepsilon}{\bar\varepsilon}(z-\hat x^j)$.
Additionally we use that $x/\varepsilon=z/\bar\varepsilon+I$ and $\phi$ is
$Y$--periodic. (\ref{P8C approximation 1}) is equal to
$$
\sum_{j=1}^N\sum_{I\in{\bf J}_\varepsilon^{int}(j)}
\left(\frac{\varepsilon}{\bar\varepsilon}\right)^d\int\limits_{W_j} 
a_M(z)\phi\Bigl(x^j_I(z),\frac{z}{\bar\varepsilon}\Bigr)\, dz=
$$
$$
=\sum_{j=1}^N\frac{1}{\mu(W_j)}\int\limits_{W_j}a_M(z)\Bigl[
\sum_{I\in{\bf J}_\varepsilon^{int}(j)}\varepsilon^d
\phi\Bigl(x^j_I(z),\frac{z}{\bar\varepsilon}\Bigr)\Bigr]\, dz
$$
and it is approximately equal to (\ref{P8C RHS}).
If $z\in W_j$ then  $x^j_I(z)$ is some point in $\varepsilon$-cube
$\Delta^j_I$. For small enough $\varepsilon$ the integral
from the continuous function $\phi$ over $\Omega_j$ can be approximated by a sum
with an error not larger than $\delta$ for all $z\in W_j$, $j\in{1,\dots,N}$:
$$
|\int_{\Omega_j}\phi\Bigl(x,\frac{z}{\bar\varepsilon}\Bigr)\, dx-
\sum_{I\in{\bf J}_\varepsilon^{int}(j)}\varepsilon^d
\phi\Bigl(x^j_I(z),\frac{z}{\bar\varepsilon}\Bigr)|\le
\delta=\frac{E\bar\varepsilon^d}{2N\|a_M\|_{L^1(\tilde\Omega)}}
$$
Finally, for small enough $\varepsilon$
$$
|\int\limits_\Omega\int\limits_Y a(x,y)\phi(x,y)\,dy\,dx-
\int\limits_{\Omega}a(x,\frac{x}{\varepsilon})\phi(x,\frac{x}{\varepsilon})\, dx|\le
\frac{E}{2}+\sum_{j=1}^N\frac{\delta}{\mu(W_j)}\int\limits_{W_j}|a_M(z)|\, dz
\le E.
$$
\end{proof}
\section{Admissibility of the ${\mathcal Continuous}$, 
${\mathcal Discrete}$ Extensions}
\label{s:Admissibility of the Continuous Discrete Extensions}
Starting from here, if it is not explicitly mentioned, then a 
two-scale extension means either the ${\mathcal Continuous}$ or the
${\mathcal Discrete}$ Extensions as defined in Subsections~\ref{Section Continuous Extension},
\ref{Section Discrete Extension}.
\begin{corollary}
\label{Cor. 1}
If $a_M(\cdot)\in L^p(\tilde\Omega)$, $p\in\mathbb{N}$, 
$\phi\in C(\overline\Omega\times\overline Y)$ is $Y$--periodic,  
$a(x,y)$ is either the ${\mathcal Continuous}$ or the ${\mathcal Discrete}$ Extension of $a_M(x)$
then 
$$
\lim\limits_{\varepsilon\to 0}\int_\Omega 
a(x,x/\varepsilon)^p\phi(x,x/\varepsilon)\, dx=
\int_\Omega\int_Y a(x,y)^p\phi(x,y)\,dx\,dy. 
$$
\end{corollary}
\begin{proof}
$|a_M|^p=|a_M^p|$ then $b_M(\cdot)=a_M(\cdot)^p\in L^1(\tilde\Omega)$.
As we know from Cor.~\ref{P0} the two-scale extension for $b_M(\cdot)$
is $b(x,y)=a(x,y)^p$. Then we can apply Lem.~\ref{P8 B} or \ref{P8 C} for $b_M(\cdot)$. 
\end{proof}
\begin{corollary}
\label{Cor. 2}
If $a_M(\cdot)\in L^2(\tilde\Omega)$, $\psi\in C(\overline\Omega\times\overline Y)$ is $Y$--periodic
then $a(x,y)\psi(x,y)$ is an "admissible" test function .
\end{corollary}
\begin{proof}
We should check (\ref{ATF condition 1}), (\ref{ATF condition 2}) for $a(x,y)\psi(x,y)$.
$\tilde\Omega$ is bounded; therefore $a_M(\cdot)\in L^1(\tilde\Omega)$.
For arbitrary $\phi\in{\cal B}_{TF}$ we can choose $p=1$ and $\phi\psi\in C(\overline\Omega\times\overline Y)$
instead of $\phi$ in Cor.~\ref{Cor. 1} to verify (\ref{ATF condition 1}).
The second condition (\ref{ATF condition 2}) is again a consequence of Cor.~\ref{Cor. 1} with $p=2$
and $\phi=\psi^2\in C(\overline\Omega\times\overline Y)$. 
\end{proof}
\section{Application to the elliptic equation}
\label{Section Application to the elliptic equation}
We return back to the practical problem from the Subsection 
~\ref{Subsection Elliptic example}.

In the context of two-scale convergence, the sequence of problems 
(\ref{sequence P(varepsilon)}) was investigated in ~\cite{Allaire},\S 2.
Now what is required is to go through the proofs in order to convince ourselves that
they still work in our case when $a(x,y)$ is a two--scale extension of $a_M(x)$. 
$$
\alpha|\xi|^2\le \xi^T a_M(x)\xi,\quad |a_M(x)\xi|\le\beta |\xi| 
\qquad\mbox{in $\tilde\Omega$, for any $\xi\in\mathbb{R}^d$}
$$
implies that (Cor.~\ref{P0})
$$
\alpha|\xi|^2\le \xi^T a(x,y)\xi,\quad |a(x,y)\xi|\le\beta |\xi|  
\qquad\mbox{in $\Omega\times Y$, for any $\xi\in\mathbb{R}^d$}
$$
For any $\varepsilon>0$, $a(x,x/\varepsilon)$ is measurable in $\Omega$
provided $a_M(\cdot)$ is measurable in $\tilde\Omega$ 
(Prop.~\ref{P4 B}, Prop.~\ref{P4 C}, $\phi\equiv 1$).
Therefore for $f\in L^2(\Omega)$ the problems (\ref{sequence P(varepsilon)})
are uniquely solvable and their solutions are uniformly bounded in $H^1_0(\Omega)$.
\begin{theorem}
The sequence $u_\varepsilon$ of  solutions of (\ref{sequence P(varepsilon)})
converges weakly in $H^1_0(\Omega)$ (and strongly in $L^2(\Omega)$) to $u_0$,
a unique solution of the limit problem:
\begin{equation}
\label{averaged elliptic problem}
\mbox{${\mathcal P}^{0}$:}\qquad
-\nabla\cdot(A(x)\nabla u_0)=f\qquad\mbox{in $\Omega$,}
\qquad u_0|_{\partial\Omega}=0.
\end{equation}
where
\begin{equation}
\label{Averaged coefficient Aij}
A_{ij}(x)=\int_Y e_i^T a(x,y)\Bigl(\nabla_y w_j(x,y)+e_j\Bigr)\,dy,
\end{equation}
$e_j$ - basis vectors, $w_j(x,y)$ ($j=1,\dots ,d$) are solutions of 
the cell problems:
\begin{equation}
\label{Cell Problem}
\left\{
\begin{array}{l}
-\nabla_y\cdot\Bigl(a(x,y)\bigl(\nabla_y w_j(x,y)+e_j\bigr)\Bigr)=0
\qquad\mbox{in Y}\\
\int\limits_Y w_j(x,y)\,dy=0,\qquad\mbox{$w_j(x,y)$ is $Y$-periodic in $y$} 
\end{array}
\right.
\end{equation}
\end{theorem}
\begin{proof}
For bounded $\tilde\Omega$ we have $L^\infty(\tilde\Omega)\subset L^2(\tilde\Omega)$.
$\phi(x)\in D(\Omega)$, $\phi_1(x,y)\in D[\Omega; C^\infty_{per}(Y)]$ 
have continuous derivatives
in $\overline\Omega\times\overline Y$ and according to Cor.~\ref{Cor. 2}
$$
[\nabla\phi(x)+\nabla_y\phi_1(x,y)]^T a(x,y),\quad
[\nabla_x\phi_1(x,y) ]^T a(x,y)
$$
are row vectors consisting of admissible test functions. Therefore
it is still possible in this case to pass to the two-scale limit in
~\cite{Allaire}, (2.10) to obtain ~\cite{Allaire}, (2.11).
The remaining part is given by ~\cite{Allaire}, Proof of T.2.3.
\\  
We note that the uniqueness of the solution to the limit problem resulting
in the convergence of the whole sequence (not just some subsequence) is important
to insure that the solution to the initial problem 
${\mathcal P}={\mathcal P}^{\textstyle\bar\varepsilon}$ belongs to the 
convergent sequence.
\end{proof}

What can we say about the averaged coefficient $A$?
In the case of the ${\mathcal Continuous}$ Extension $a(x,y)$ this coefficient
should be calculated at each point $x\in\Omega$ and it depends on the initial coefficient
$a_M(\cdot)$ in $\bar\varepsilon$-cube $W(x)$ around $x$.
What happens with $A$ if we slightly move from $x$ to $x+h$?
For small enough $h$ the volume $W(x)$ has a large intersection with $W(x+h)$ and 
consequently the coefficients $a(x,\cdot)$, $a(x+h,\cdot)$, 
which play a crucial role in the cell problem,
differ from each other only in a small volume. 
Therefore continuity of the averaged coefficient depends on the form
of cell problem.
\begin{proposition}
The coefficient $A(x)$ calculated from the ${\mathcal Continuous}$ Extension $a(x,y)$
is continuous in $\overline\Omega$.
\end{proposition}
\begin{proof}
Although our cell problem (needed for calculation of $A(x)$) 
is formulated in $Y$ and has a variational form: find $w_j(x,\cdot)\in H^1_{per}(Y)\setminus\mathbb{R}$ such that
$$
\int_Y \nabla_y\phi(y)^Ta(x,y)\nabla_y w_j(x,y)\, dy=
-\int_Y\nabla_y\phi(y)^Ta(x,y)e_j\, dy\qquad
\forall\phi\in H^1_{per}(Y)\setminus\mathbb{R},
$$
we prefer to deal with the cell problem in terms of $a_M(\cdot)$, not in 
$a(\cdot,\cdot)$. To do this it is better to substitute 
$\bar\varepsilon$-periodic functions for $Y$-periodic:
$$
\tilde\phi(z):=\phi(\frac{z}{\bar\varepsilon}),\quad
\tilde w_j(x,z):=w_j(x,\frac{z}{\bar\varepsilon}),\quad
\tilde a(x,z)=a(x,\frac{z}{\bar\varepsilon})
$$
After this substitution the integrals will be over $\bar\varepsilon Y$ from
$\bar\varepsilon$-periodic functions in $z$. Therefore they are equal to integrals
over $W(x)$, where $\tilde a(x,z)=a_M(z)$.
 
In new terms the problem has the form: 
find $\tilde w_j(x,\cdot)\in H^1_{per}(W(x))\setminus\mathbb{R}$ such that
for all $\tilde\phi(z)\in H^1_{per}(W(x))\setminus\mathbb{R}$ holds the equality\quad
${\mathcal B}(\tilde w_j(x,\cdot),\tilde\phi)={\mathcal L}_j(\tilde\phi)$,\quad where
\begin{equation}
\label{Cell problem A(x) var. form} 
{\mathcal B}(\tilde w,\tilde\phi):=
\int_{W(x)}\nabla_z\tilde\phi(z)^T a_M(z)\nabla_z\tilde w(z)\, dz,\quad
{\mathcal L}_j(\tilde\phi):=
-\frac{1}{\bar\varepsilon}\int_{W(x)}\nabla_z\tilde\phi(z)^T a_M(z)e_j\, dz. 
\end{equation}
The $\bar\varepsilon$-periodic function from $H^1_{per}(W(x))\setminus\mathbb{R}$
is also a function from $H^1_{per}(W)\setminus\mathbb{R}$, where $W$ is an arbitrary 
$\bar\varepsilon$-cube, 
$\|\phi\|_{H^1_{per}(W(x))\setminus\mathbb{R}}=\|\phi\|_{H^1_{per}(W)\setminus\mathbb{R}}$.
The Poincare inequality:
for all $\phi\in H^1_{per}(Y)\setminus\mathbb{R}$ and 
$\tilde\phi\in H^1_{per}(W)\setminus\mathbb{R}$ 
$$
\int_Y\phi(y)^2\,dy\le C^2_\#\int_Y|\nabla_y\phi(y)|^2\, dy,\qquad
\int_W\tilde\phi(z)^2\,dz\le 
C^2_\#\bar\varepsilon^2\int_W|\nabla_z\tilde\phi(z)|^2\, dz.
$$
We note that the 'small' size of $\bar\varepsilon$ is not important here. 
It is just a fixed constant.
The bilinear form ${\mathcal B}$ is elliptic and bounded on 
$H^1_{per}(W)\setminus\mathbb{R}$:
$$
{\mathcal B}(\tilde w,\tilde w)\ge\alpha\int_{W(x)}|\nabla_z\tilde w(z)|^2\,dz\ge
\frac{\alpha}{1+C^2_\#\bar\varepsilon^2}\|\tilde w\|^2_{H^1_{per}(W)\setminus\mathbb{R}},
$$
$$
|{\mathcal B}(\tilde w,\tilde\phi)|\le
\beta 
\|\tilde w\|_{H^1_{per}(W)\setminus\mathbb{R}}
\|\tilde\phi\|_{H^1_{per}(W)\setminus\mathbb{R}}.
$$
The linear functional ${\mathcal L}_j$ is bounded on $H^1_{per}(W)\setminus\mathbb{R}$:
$$
|{\mathcal L}_j(\tilde\phi)|\le\beta\bar\varepsilon^{d/2-1}
\|\tilde\phi\|_{H^1_{per}(W)\setminus\mathbb{R}}.
$$
Therefore the cell problem has a unique solution 
$\tilde w_j(x,\cdot)\in H^1_{per}(W)\setminus\mathbb{R}$, satisfying:
$$
\|\tilde w_j(x,\cdot)\|_{H^1_{per}(W)\setminus\mathbb{R}}\le
\beta\bar\varepsilon^{d/2-1}\frac{1+C^2_\#\bar\varepsilon^2}{\alpha} 
$$
The formula (\ref{Averaged coefficient Aij}) written in terms of $\tilde w_j$ is
\begin{equation}
\label{Averaged coefficient Aij from W}
A_{ij}(x)=\frac{1}{\bar\varepsilon^d}\int_{W(x)}e_i^T a_M(z)
(\bar\varepsilon\nabla_z\tilde w_j(x,z)+e_j)\, dz.
\end{equation}
$x$ is an arbitrary point from $\overline\Omega$. To check the continuity of
$A(\cdot)$ we fix some point $x\in\overline\Omega$ and consider some point 
$x+h\in\overline\Omega$ to compare $A(x)$ and $A(x+h)$. 
The cell problem for $A(x+h)$ is
\begin{equation}
\label{Cell problem A(x+h) var. form}
\int_{W(x+h)} 
\nabla_z\tilde\phi(z)^T a_M(z)\nabla_z\tilde w_j(x+h,z)\, dz=
-\frac{1}{\bar\varepsilon}\int_{W(x+h)}\nabla_z\tilde\phi(z)^T a_M(z)e_j\, dz 
\end{equation}
$S_+:=W(x)\setminus W(x+h)$, $S_-:=W(x+h)\setminus W(x)$. 
$W(x)=W(x+h)\cup S_+\setminus S_-$. 
$\int_{W(x)}=\int_{W(x+h)}+\int_{S_+}-\int_{S_-}$. The problem 
(\ref{Cell problem A(x) var. form}) can be re-written:
$$
\int_{W(x+h)}\nabla_z\tilde\phi(z)^T a_M(z)\nabla_z\tilde w_j(x,z)\, dz=
-\frac{1}{\bar\varepsilon}\int_{W(x)}\nabla_z\tilde\phi(z)^T a_M(z)e_j\, dz+
$$
\begin{equation}
\label{Cell problem A(x) var. form similar to A(x+h)}
+\int_{S_-}\nabla_z\tilde\phi(z)^T a_M(z)\nabla_z\tilde w_j(x,z)\, dz
-\int_{S_+}\nabla_z\tilde\phi(z)^T a_M(z)\nabla_z\tilde w_j(x,z)\, dz
\end{equation}
We substitute (\ref{Cell problem A(x+h) var. form}) from 
(\ref{Cell problem A(x) var. form similar to A(x+h)}) denoting 
$\theta(z):=\tilde w_j(x,z)-\tilde w_j(x+h,z)$:
\begin{equation}
\label{problem for theta var. form}
\int_{W(x+h)}\nabla_z\tilde\phi(z)^T a_M(z)\nabla_z\theta(z)\, dz=
\hat{\mathcal L}(\phi),
\end{equation}
$$
\hat{\mathcal L}(\phi)=
\frac{1}{\bar\varepsilon}\int_{S_-}\nabla_z\tilde\phi(z)^T a_M(z)e_j\, dz
-\frac{1}{\bar\varepsilon}\int_{S_+}\nabla_z\tilde\phi(z)^T a_M(z)e_j\, dz+
$$
$$
+\int_{S_-}\nabla_z\tilde\phi(z)^T a_M(z)\nabla_z\tilde w_j(x,z)\, dz
-\int_{S_+}\nabla_z\tilde\phi(z)^T a_M(z)\nabla_z\tilde w_j(x,z)\, dz
$$
The point $x$ was fixed. Consequently the function $\tilde w_j(x,\cdot)$ is also a fixed
function. $\theta$ belongs to $H^1_{per}(W)\setminus\mathbb{R}$.
Its norm can be estimated by treating (\ref{problem for theta var. form}) as a 
variational problem for the unknown $\theta$.
$$
\Bigl|
\frac{1}{\bar\varepsilon}\int_{S_+}\nabla_z\tilde\phi(z)^T a_M(z)e_j\, dz\Bigr|
\le\beta 
\|\tilde\phi\|_{H^1_{per}(W)\setminus\mathbb{R}}\frac{\sqrt{\mu(S_+)}}{\bar\varepsilon}
$$
$$
\Bigl|
\int_{S_+}\nabla_z\tilde\phi(z)^T a_M(z)\nabla_z\tilde w_j(x,z)\, dz\Bigr|\le
\beta 
\int_{S_+}|\nabla_z\tilde\phi(z)||\nabla_z\tilde w_j(x,z)|\, dz\le
$$
$$
\le\beta 
\|\nabla_z\tilde\phi(z)\|_{[L^2(S_+)]^d}
\sqrt{\int_{S_+}|\nabla_z\tilde w_j(x,z)|^2\, dz}
$$
$|\nabla_z\tilde w_j(x,\cdot)|$ is a fixed function from $L^2(W(x))$.
Due to the absolute continuity of the Lebesgue integral for arbitrary
$E_h>0$ one can find such $\delta>0$ that the integral from 
$|\nabla_z\tilde w_j(x,\cdot)|^2$ over any set in $W(x)$ 
is less than $E_h$ if the set's measure is less than $\delta$. 
The measure of $S_+$, is arbitrarily small provided $h$ is small enough.
$$
\Bigl|
\int_{S_+}\nabla_z\tilde\phi(z)^T a_M(z)\nabla_z\tilde w_j(x,z)\, dz\Bigr|
\le\beta 
\|\tilde\phi\|_{H^1_{per}(W)\setminus\mathbb{R}}\sqrt{E_h}
$$
Similar estimations can be done for the integrals over $S_-$.
Therefore $\|\hat{\mathcal L}\|$ and consequently $\|\theta\|_{H^1_{per}(W)\setminus\mathbb{R}}$
are arbitrarily small for small enough $h$.
$$
\|\theta\|_{H^1_{per}(W)\setminus\mathbb{R}}\le
\frac{1+C_\#^2\bar\varepsilon^2}{\alpha}\|\hat{\mathcal L}\|.
$$
This helps us to estimate $|A_{ij}(x)-A_{ij}(x+h)|$.
$$
A_{ij}(x)-A_{ij}(x+h)=\frac{1}{\bar\varepsilon^d}\int_{W(x+h)}e_i^T a_M(z)
\bar\varepsilon\nabla_z\theta(z)\, dz+
$$
$$
+\frac{1}{\bar\varepsilon^d}\int_{S_+}e_i^T a_M(z)
(\bar\varepsilon\nabla_z\tilde w_j(x,z)+e_j)\, dz
-\frac{1}{\bar\varepsilon^d}\int_{S_-}e_i^T a_M(z)
(\bar\varepsilon\nabla_z\tilde w_j(x,z)+e_j)\, dz
$$
The absolute value of the first term and the integral over $S_+$ splitted 
into two parts
$$
|\frac{1}{\bar\varepsilon^d}\int_{W(x+h)}e_i^T a_M(z)\bar\varepsilon\nabla_z\theta(z)\, dz|\le
\left.\beta 
\|\theta\|_{H^1_{per}(W)\setminus\mathbb{R}}\right/\bar\varepsilon^{d/2-1},
$$
$$
|\frac{1}{\bar\varepsilon^d}\int_{S_+}e_i^T a_M(z)
\bar\varepsilon\nabla_z\tilde w_j(x,z)\, dz|\le
\left.\beta 
\sqrt{\mu(S_+)E_h}\right/\bar\varepsilon^{d-1},
$$
$$
|\frac{1}{\bar\varepsilon^d}\int_{S_+}e_i^T a_M(z)e_j\, dz|\le\left.
\beta 
\mu(S_+)\right/\bar\varepsilon^d
$$
can be made arbitrarily small (due to the small terms 
$\|\theta\|_{H^1_{per}(W)\setminus\mathbb{R}}$, $E_h$, $\mu(S_+)$)
by choosing small enough $h$. Together with similar estimations for $S_-$ we have 
the continuity of $A_{ij}(x)$. 
\end{proof}

In the case of the ${\mathcal Discrete}$ Extension $a(x,y)$, 
the averaged coefficient $A(x)$ is constant in each $\Omega_k$ ($k=1,\dots, N$). 
To determine it one has to solve $N$ cell problems 
(\ref{Cell problem A(x) var. form}),(\ref{Averaged coefficient Aij from W})
with $x=\hat x^k$ (centers of $W_k$).
This case is realizable in comparison to solving the cell problems at each point in 
$\Omega$. On the other hand the averaged coefficient being continuous can be interpolated between a finite number of points where it is calculated via cell problems.
Here one should be careful since for small $\bar\varepsilon$ the averaged coefficient $A(x)$ is as oscillatory as the initial coefficient $a_M(\cdot)$. Increasing 
$\bar\varepsilon$ we expect $A(x)$ to become a function with more and more slow variations and in the subdomains of $\Omega$ where the coefficient 
$a_M(\cdot)$ can be classified as 'spatially homogeneous' it might be close to a constant coefficient.

In Section.~\ref{Section Three approaches to construct a two--scale extension}
we had a restriction on $\bar\varepsilon$ from above: $\bar\varepsilon$ should be small in
comparison with the typical size of $\Omega$. Solving the limit problem numerically with
some typical discretization step $h$ provides a restriction for $\bar\varepsilon$ 
from below: roughly speaking, $\bar\varepsilon$ should not be smaller than $h$.

Solving numerically the large number of cell problems is a time consuming task, which
can be done in parallel since cell problems are independent from each
other and the limit problem.
The computational resources can be also saved at least in the following cases:\\  
$\bullet$ $a_M(x)$ has slow variations (for example it can be a constant) in some subdomain 
$\Omega_{sv}\subset\Omega$. Then inside $\Omega_{sv}$ there is no need
to average.\\
$\bullet$ $a_M(x)$ is $\bar\varepsilon$-periodic in $\Omega_\#\subset\Omega$ and 
the directions of periodicity coincide with coordinate axes. Then the constant averaged 
coefficient inside $\Omega_\#$ can be calculated by solving only one cell problem.

Additionally one can also try to combine this with other types of averaging:\\
$\bullet$ if the micro coefficient can be classified as statistically homogeneous in 
some subdomain $\Omega_{sh}$ with known averaged value $A_{sh}$ or\\
$\bullet$ if the averaged coefficient $A_{ed}$ in $\Omega_{ed}$ is experimentally 
determined.\\
In these cases one can use the coefficients $A_{sh}$ inside $\Omega_{sh}$ and
$A_{ed}$ in $\Omega_{ed}$ instead of solving cell-problems there. 
\section{Some concluding remarks}
The coefficient $a(x,x/\varepsilon)$ is often used in homogenization as a
generalization of the periodic coefficient $a(x/\varepsilon)$.
In this paper we propose a way to correspond an averaged (limit) problem for
the initial microscopical problem with non-periodic rapidly oscillated
coefficient using the results from homogenization
together with a special choice of $a(x,y)$ -- the two-scale extension of the initial
coefficient. The results from homogenization (if not formal) usually require
some conditions (like its smoothness) on $a(x,y)$. The lack of smoothness e.g. in the two-scale convergence method can be partially compensated 
by the ''admissibility'' of the two-scale extensions, so that e.g. for the 
second order elliptic equation the convergence of $u_\varepsilon$ to $u_0$
still holds as in the periodic case.

\begin{figure}[tp]
\begin{center}
\scalebox{0.9}{\includegraphics{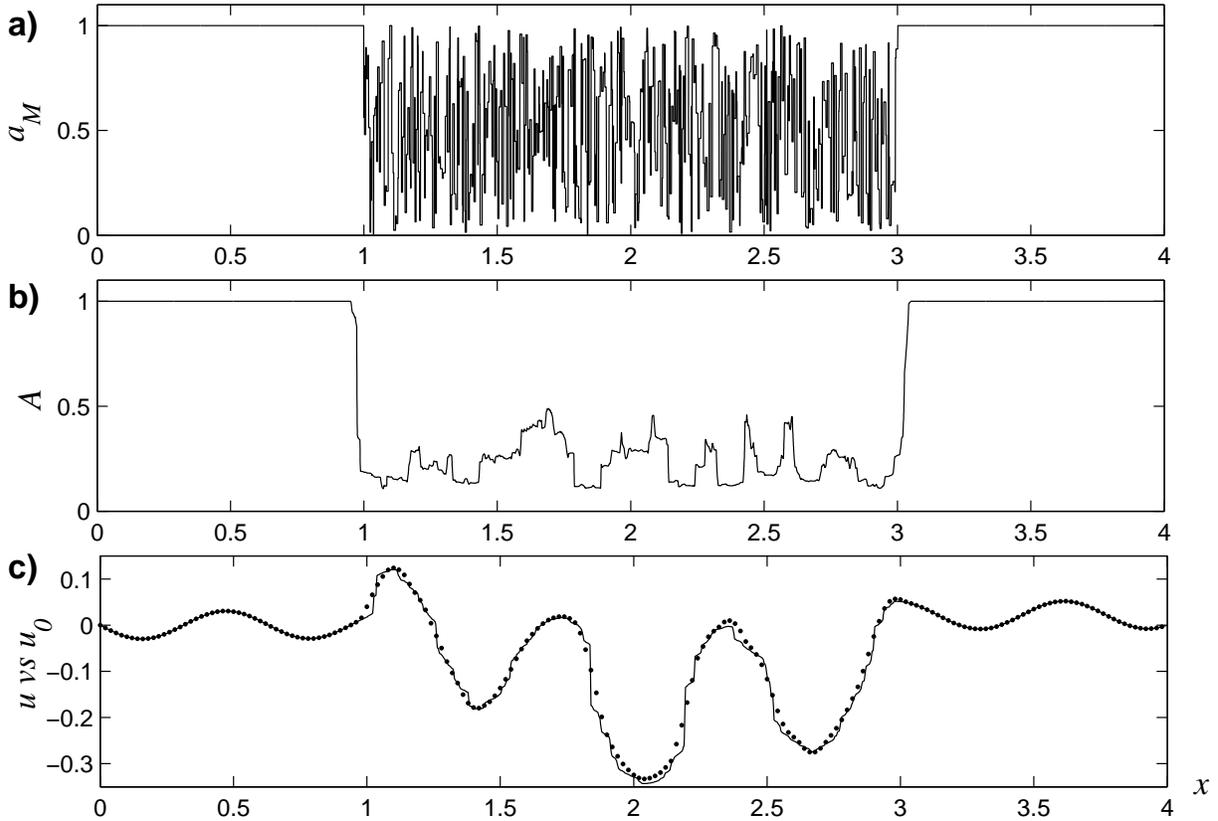}}
\end{center}
\caption{a) initial coefficient $a_M(\cdot)$; b) averaged coefficient $A(\cdot)$; c) comparison $u(x)$ with $u_0(x)$}
\label{numerical example}
\end{figure}

To show that this approach can be useful we present here a 1D example where 
$a_M(x)$ from (\ref{initial elliptic problem}) and $A(x)$ 
from (\ref{averaged elliptic problem}) are plotted in 
Fig.~\ref{numerical example}(a,b) respectively. 
To calculate $A(\cdot)$, the ${\mathcal Continuous}$
extension for $\bar\varepsilon=0.1$ was used.
The semi-analytical solutions $u(x)$ (solid line) and $u_0(x)$ (dots) corresponding to $f(x)=-3\sin(10x)$ are compared in 
Fig.~\ref{numerical example}(c). 

In a 2D test presented in
\cite{Laptev GAMM} a fine scale reference solution to ${\mathcal P}$ is compared
with a $H^1$-corrected coarse solution to ${\mathcal P}^{0}$ (this
classical correction is described e.g. in \cite[p.76]{BLP}).
$A(x)$ in ${\mathcal P}^{0}$ is calculated via the ${\mathcal
  Discrete}$ extension of a randomly generated smooth 2D function $a_M(\cdot)$.

In consequent publications we are planning to present numerical
results for the elliptic problem in 1D and 2D more systematically
together with some other two-scale extensions.

\end{document}